\documentclass[12pt]{amsart}

\newtheorem{theo}{Theorem}[section]

\newtheorem{prop}{Proposition}[section]

\usepackage{amsfonts}
\usepackage{amssymb}
\usepackage{eucal}
\numberwithin{equation}{section}
\title
[CONDITIONS DETERMINED BY THE OPOZDA-VERSTRAELEN AFFINE CURVATURE TENSOR]{ON 
HYPERSURFACES SATISFYING CONDITIONS DETERMINED 
BY THE OPOZDA-VERSTRAELEN AFFINE CURVATURE TENSOR}

\author{Ryszard Deszcz, Ma\l gorzata G\l ogowska and Marian Hotlo\'{s}}
\begin{document}

\dedicatory{Dedicated to Professor Udo Simon on his eighty-second birthday}

\begin{abstract}
Using the Blaschke-Berwald metric and the affine shape operator of a hypersurface $M$ 
in the $(n+1)$-dimensional real affine space
we can define some generalized curvature tensor named the Opozda-Verstraelen affine curvature tensor. 
In this paper we determine curvature conditions of pseudosymmetry type expressed by this tensor 
for locally strongly convex hypersurfaces $M$, $n > 2$, with two distinct affine principal curvatures or
with three distinct affine principal curvatures assuming 
that at least one affine principal curvature has multiplicity $1$.
\end{abstract}

\maketitle

\section{Introduction}

Let $\mathbb{A}^{n+1}$, $n \geq 2$, be the standard $(n+1)$-dimensional real affine space, 
i.e. $\mathbb{R}^{n+1}$ equipped with its standard flat connection $\widetilde{\nabla}$
and the volume element $\widetilde{\Theta}$ given by the determinant.
Since $\widetilde{\nabla}$ has no torsion and $\widetilde{\Theta}$ is parallel with respect to $\widetilde{\nabla}$,
the pair $(\widetilde{\nabla},\widetilde{\Theta})$ determines an equiaffine structure on $\mathbb{R}^{n+1}$.
The space $\mathbb{A}^{n+1}$ is a homogeneous space under the natural action of the unimodular affine group 
$SL(n+1, \mathbb{R}) \propto \mathbb{R}^{n+1}$ {\cite[Preliminaries] {1990_OV1}}.
Further, let $M$ be a connected non-degenerate hypersurface in $\mathbb{A}^{n+1}$ 
with the affine normal $\xi$
and the induced equiaffine structure $( \nabla, \theta )$ 
(see, e.g., \cite{{LSZ}, {LSZH}, {Nomizu-Sasaki}, {1990_OV1}, {1985_VV}, {1984_V}}).
Using the curvature tensor $\mathcal{R}$ of the connection $\nabla$, 
the affine shape operator $\mathcal{S}$, the Blaschke-Berwald metric $h$ of $M$,
and the Gauss and the Ricci equations of $M$ in $\mathbb{A}^{n+1}$,
\begin{eqnarray*}
\mathcal{R}(X,Y)Z\ =\  h(Y,Z) \mathcal{S} X - h(X,Z) \mathcal{S} Y\  \  \mbox{and}\  \  h(X, \mathcal{S}Y)\  =\  h(Y, \mathcal{S}X), 
\end{eqnarray*}
we can define
the ge\-ne\-ralized curvature tensor $R^{\ast}$  
named the {\sl{Opozda-Verstraelen affine curvature tensor}}  of $M$
{\cite{{1990_OV1}, {1999_OV2}}}. Namely, we have 
\begin{eqnarray}
R^{\ast} (X_{1}, \ldots , X_{4}) &=& h( \mathcal{R}( X_{1}, X_{2})\mathcal{S} X_{3} , X_{4})
\ =\
h ( h(X_{2}, \mathcal{S} X_{3}) \mathcal{S} X_{1} - h(X_{1}, \mathcal{S} X_{3}) \mathcal{S} X_{2},  X_{4})\nonumber\\
&=&
h(X_{2}, \mathcal{S} X_{3}) h ( \mathcal{S}X_{1} , X_{4})
-
h(X_{1}, \mathcal{S} X_{3})
h( \mathcal{S}X_{2} , X_{4})\nonumber\\
&=&
h(X_{1}, \mathcal{S}X_{4} )
h(X_{2}, \mathcal{S}X_{3} ) 
-
h(X_{1}, \mathcal{S}X_{3} )
h(X_{2}, \mathcal{S}X_{4} ) ,
\label{aff03aa}
\end{eqnarray}

\hspace*{0pt}\hrulefill\hspace*{155mm}\\
{\footnotesize{\indent
{\bf{Mathematics Subject Classification (2010):}} 
53A15, 53B20, 53B25.

{\bf{Key words and phrases:}} partially Einstein space, pseudosymmetry type curvature condition, Roter space,
real affine space, non-degenerate hypersurface, Blaschke-Berwald metric, 
affine quasi-umbilical hypersurface, affine $2$-quasi-umbilical hypersurface, 
affine principal curvature,
Opozda-Verstraelen affine curvature tensor,
affine quasi-umbilical$^{\ast}$ hypersurface, affine $2$-quasi-umbilical$^{\ast}$ hypersurface, 
affine Einstein$^{\ast}$ hypersurface,
affine partially Einstein$^{\ast}$ hypersurface, hypersurface, $2$-quasi-umbilical hypersurface,
principal curvature.}}

\newpage

\noindent
where $X,Y,Z, X_{1}, \ldots , X_{4}$ are tangent vector fields on $M$. 
The tensor $R^{\ast}$ will be also called the {\sl{Opozda-Verstraelen tensor}}. If we set
\begin{eqnarray}  
S(X,Y) &=& h(X,\mathcal{S} Y)  
\label{aff03zzz} 
\end{eqnarray}
then (\ref{aff03aa}) turns into 
\begin{eqnarray*}
R^{\ast} (X_{1}, \ldots , X_{4}) &=& S(X_{1},X_{4}) S(X_{2},X_{3}) - S(X_{1},X_{3}) S(X_{2},X_{4}) ,
\end{eqnarray*}
i.e., in short, we have
%\noindent
\begin{eqnarray}
R^{\ast} &=& \frac{1}{2}\, S \wedge S . 
\label{aff03} 
\end{eqnarray}

We refer to sections 2 and 3 for precise definitions of the symbols used. 
These sections also contain preliminary results. 
We mention that in sections 2, 4 and 5
we will denote by $S$ the Ricci tensor of the considered semi-Riemannian manifolds $(M,g)$. 
In Section 4  we present a survey on semi-Riemannian manifolds satisfying curvature conditions named conditions of pseudosymmetry type.
For a more wider presentation on such conditions we refer to {\cite[Sections 1 and 3] {2016_DGJZ}} (see also {\cite[Section 1] {2018_DGZ}}).
These conditions determine certain classes of semi-Riemannian manifolds, for instance: 
pseudosymmetric manifolds, Ricci-pseudosymmetric manifolds and manifolds with pseudosymmetric Weyl tensor.
Curvature tensor $R$ of some semi-Riemannian manifolds $(M,g)$ 
is a linear combination of the Kulkarni-Nomizu products formed by Ricci tensor $S$ and the metric $g$,
as $S \wedge S$, $g \wedge S$ and $g \wedge g$.
Such manifolds are called Roter type manifolds, or Roter manifolds, or Roter spaces.     
Section 5 is related to that class of manifolds. From Proposition 3.2 it follows  
that on any Roter space various conditions of pseudosymmetry type are satisfied.      
For a more wider presentation on Roter type manifolds we refer to {\cite[Sections 1 and 3] {2016_DGJZ}}.

Let $\mathrm{Ric} (R^{\ast})$, $\kappa (R^{\ast})$ and 
\begin{eqnarray}
\mathrm{Weyl} (R^{\ast}) &=& R^{\ast} - \frac{1}{n-2}\, h \wedge \mathrm{Ric}(R^{\ast}) 
+ \frac{ \kappa ( R^{\ast} ) }{2 (n-2) (n-1)}\, h \wedge h ,
\label{aff55aff55}
\end{eqnarray}
be the Ricci tensor, the scalar curvature and the Weyl tensor 
determined by the metric $h$ and the tensor $R^{\ast}$, respectively \cite{1990_OV1}. 
From (\ref{aff03}), by suitable contractions with respect to $h$, we get (cf. {\cite[eq. (2.12)] {1990_OV1}})
\begin{eqnarray}
(a)\ \ \mathrm{Ric}(R^{\ast})
 \ =\ \mathrm{tr}_{h}(\mathcal{S})\, S - S^{2} ,
\ \ & &\ \
(b)\ \ \kappa ( R^{\ast} ) \ =\ (\mathrm{tr}_{h}(\mathcal{S}))^{2} - \mathrm{tr}_{h}(\mathcal{S}^{2}) ,
\label{aff0404} 
\end{eqnarray}
where $h_{ij}$, $h^{ij}$, $S_{ij}$ and $S^{2}_{ij} = S_{ik}h^{kl}S_{lj}$
are the local components of the tensors $h$, $h^{-1}$, $S$ and $S^{2}$, respectively.
It is obvious that the tensors $R^{\ast}$ and $\mathrm{Weyl} (R^{\ast})$ are generalized curvature tensors. 
For these tensors and the tensors $h$ and 
$\mathrm{Ric}(R^{\ast})$ we can define the following $(0,6)$-tensors
(see sections 2 and 3):
$R^{\ast} \cdot R^{\ast}$,
$R^{\ast} \cdot \mathrm{Weyl} (R^{\ast})$,
$\mathrm{Weyl} (R^{\ast}) \cdot R^{\ast}$,
$\mathrm{Weyl} (R^{\ast}) \cdot \mathrm{Weyl} (R^{\ast})$,
$Q( h , R^{\ast})$,
$Q(\mathrm{Ric} (R^{\ast}), R^{\ast})$,
$Q( h , \mathrm{Weyl} (R^{\ast}) )$ and
$Q( \mathrm{Ric} (R^{\ast}) , \mathrm{Weyl} (R^{\ast}))$.
Non-degenerate hypersurfaces $M$ in $\mathbb{A}^{n+1}$, $n \geq 3$,
satisfying conditions imposed on $R^{\ast}$ and some tensors obtained from $R^{\ast}$
were investigated in \cite{{1991D1}, {1992D2}, {1990_OV1}, {1999_OV2}}. 
For instance, in \cite{1992D2} (see also Theorem 6.1)
it was proved that the following identity is satisfied on any 
hypersurface $M$ in $\mathbb{A}^{n+1}$, $n \geq 3$, 
\begin{eqnarray}
R^{\ast} \cdot R^{\ast} &=& Q(\mathrm{Ric} (R^{\ast}), R^{\ast}) . 
\label{aff05}
\end{eqnarray}

We denote by $\mathcal{U}_{\mathcal{S}}$ the set of all points $x$ 
of a non-degenerate hypersurface $M$ in $\mathbb{A}^{n+1}$, $n \geq 3$,
at which the affine shape operator $\mathcal{S}_{x}$
is not proportional to the identity transformation $I\!d_{x}$ at this point. 
The hypersurface $M$ in $\mathbb{A}^{n+1}$, $n \geq 3$,
is said to be {\sl affine quasi-umbilical} 
at $x \in \mathcal{U}_{\mathcal{S}}$ if at this point \cite{1989_BO1} 
\begin{eqnarray}
\mathrm{rank} (S - \rho \, h) &=& 1,\ \ \mbox{for some}\ \ \rho \in \mathbb{R} .  
\label{aff05aff}
\end{eqnarray}
For a locally strongly convex hypersurface $M$ in ${\mathbb{A}}^{n+1}$, $n \geq 2$,
$h$ is positive definite. Further, let $\rho_{1}, \rho_{2}, \ldots , \rho_{n}$ be eigenvalues 
of the affine shape operator $\mathcal{S}_{x}$ at a point $x \in M$; 
they are called the {\sl{affine principal curvatures}} (see, e.g., 
{\cite[p. 51] {LSZ}}, {\cite[p. 55] {LSZH}}).
Without loss of generality we can write that if a locally strongly convex hypersurface $M$ in ${\mathbb{A}}^{n+1}$, $n \geq 3$,
is affine quasi-umbilical at $x$ then at this point we have 
$\lambda _{1} = \rho _{1}$ and $\lambda _{2} = \rho _{2} = \cdots = \rho _{n}$.
We refer to \cite{{Dillen-Vrancken-1993}, {Dillen-Vrancken-1994}, {2014_HuLiZhang}, {2012_AJV}, {1991_Vrancken}}
for examples of affine quasi-umbilical hypersurfaces. It easy to check that 
if (\ref{aff05aff}) holds at a point $x \in \mathcal{U}_{\mathcal{S}}$ then 
(see, Section 5, eqs. (5.1) and (5.2)) 
\begin{eqnarray}
S^{2} + ( (n-2) \rho  - \mathrm{tr}_{h}(\mathcal{S})  )\, S + \rho ( \mathrm{tr}_{h}(\mathcal{S}) - (n-1) \rho )\, h &=& 0
\label{aff06aff}
\end{eqnarray}
at this point.
If $M$ is an affine quasi-umbilical hypersurface in $\mathbb{A}^{n+1}$, $n \geq 3$,
then the tensor $\mathrm{Weyl} (R^{\ast})$ vanishes \cite{1990_OV1}. 
In \cite{1992D2} 
(see also Theorem 6.2)
it was proved that the converse statement is also true, 
provided that $n \geq 4$. 
If the condition
\begin{eqnarray}
S^{2} + L_{1}\, S + L\, h &=& 0 
\label{aff51}
\end{eqnarray}
is satisfied on $M$ in $\mathbb{A}^{n+1}$, $n \geq 3$,
for some functions $L$ and $L_{1}$, then on $M$ we have (see Theorem 6.3)
\begin{eqnarray}
R^{\ast} \cdot R^{\ast} &=& L\, Q( h, R^{\ast}) .
\label{aff06}
\end{eqnarray}
Evidently, (\ref{aff06aff}) is a particular case of (\ref{aff51}).
Thus on every affine quasi-umbilical hypersurface $M$ in $\mathbb{A}^{n+1}$, $n \geq 3$, (\ref{aff06})
is satisfied (see Theorem 6.4).

According to {\cite[Section 2] {1990_OV1}}, a non-degenerate hypersurface $M$ in ${\mathbb{A}}^{n+1}$, $n \geq 2$, 
is said to be an 
{\sl{affine Einstein}}$^{\ast}$ {\sl{hypersurface}}  
if 
\begin{eqnarray} 
\mathrm{Ric} (R^{\ast}) &=& \frac{ \kappa ( R^{\ast} ) }{n}\, h  
\label{aff52} 
\end{eqnarray}
holds on $M$. In Theorem 6.6 we present curvature properties of affine Einstein$^{\ast}$ hypersurfaces.

A hypersurface $M$ in ${\mathbb{A}}^{n+1}$, $n \geq 3$, 
is called {\sl{affine partially Einstein}}$^{\ast}$ 
at a point $x \in M$ if at this point we have
\begin{eqnarray} 
(\mathrm{Ric} (R^{\ast}))^{2}
&=& 
\rho _{1} \, \mathrm{Ric} (R^{\ast}) + \rho _{2} \, h , \ \ \mbox{for some}\ \ \rho _{1}, \rho _{2} \in \mathbb{R }.  
\label{affine003} 
\end{eqnarray}
We note that from (\ref{affine003}), by (\ref{aff0404})(a), we get easily
\begin{eqnarray*}
S^{4} &=& 2 \mathrm{tr}_{h}(\mathcal{S})\, S^{3} - (\rho _{1}\, + (\mathrm{tr}_{h}(\mathcal{S}))^{2})\, S^{2} 
+ \rho _{1} \mathrm{tr}_{h}(\mathcal{S})\, S  + \rho _{2}\, h . 
%\label{affine004} 
\end{eqnarray*}

Let $M$ be a non-degenerate hypersurface $M$ in $\mathbb{A}^{n+1}$, $n \geq 4$.
Further, 
let $\mathcal{U}$ be the set of all points 
of ${\mathcal{U}}_{ \mathrm{Ric}(R^{\ast})} \cap {\mathcal{U}}_{ \mathrm{Weyl}(R^{\ast})} \subset M$ 
at which 
\begin{eqnarray}
\mathrm{rank} (\mathrm{Ric} (R^{\ast}) - \rho \, h)  \geq  2 ,\ \ \mbox{for any}\ \ \rho \in \mathbb{R} . 
\label{aff07}
\end{eqnarray}
In addition, let (\ref{aff06}) be satisfied on $\mathcal{U}$.
Using (\ref{aff05}), (\ref{aff06}) and (\ref{aff07}) we can prove that 
at every point of $\mathcal{U}$ the tensor $R^{\ast}$ is a linear combination 
of the Kulkarni-Nomizu products $\mathrm{Ric} (R^{\ast}) \wedge \mathrm{Ric} (R^{\ast})$, $h \wedge \mathrm{Ric} (R^{\ast})$
and $h \wedge h$ and the tensor 
$\mathrm{Weyl} (R^{\ast}) \cdot R^{\ast} - R^{\ast}  \cdot  \mathrm{Weyl} (R^{\ast})$
is a linear combination of the tensors  
$Q(\mathrm{Ric} (R^{\ast}), \mathrm{Weyl} (R^{\ast}) )$ 
and $Q(h, \mathrm{Weyl} (R^{\ast}))$ (see Theorem 6.5). 
Precisely, we have on $\mathcal{U}$
\begin{eqnarray}
\ \ \ \ \ \ \ 
\mathrm{Weyl} (R^{\ast}) \cdot R^{\ast} - R^{\ast}  \cdot  \mathrm{Weyl} (R^{\ast})
&=& 
Q(\mathrm{Ric} (R^{\ast}), \mathrm{Weyl} (R^{\ast}) )
- 
\frac{  \kappa (R^{\ast}) }{ n-1 }\, Q(h, \mathrm{Weyl} (R^{\ast})) .
\label{aff08}
\end{eqnarray} 

A non-degenerate hypersurface $M$ in $\mathbb{A}^{n+1}$, $n \geq 4$, 
is said to be {\sl{affine $2$-quasi-umbilical}}
at $x \in \mathcal{U}_{\mathcal{S}}$ if at this point 
\begin{eqnarray}
\mathrm{rank} (S - \rho \, h) &=& 2 ,\ \ \mbox{for some}\ \ \rho \in \mathbb{R}. 
\label{aff101}
\end{eqnarray}
Without loss of generality we can write that if a locally strongly convex hypersurface $M$ in ${\mathbb{A}}^{n+1}$, $n \geq 4$,
is affine $2$-quasi-umbilical at $x$ then at this point we have 
$\lambda _{1} = \rho _{1}$, $\lambda _{2} = \rho _{2}$ and $\lambda _{3} = \rho _{3} = \cdots = \rho _{n}$.
Theorems 7.1 and 7.2 contain results on
affine $2$-quasi-umbilical hypersurfaces. For instance, 
if some additional conditions are satisfied then the tensor $R^{\ast}$ of affine $2$-quasi-umbilical hypersurfaces 
is a linear combination of Kulkarni-Nomizu products formed by the tensors $h$, 
$\mathrm{Ric} (R^{\ast})$ and $(\mathrm{Ric} (R^{\ast}))^{2}$ (see Theorem 7.1 for details).
In Section 7 we also consider a class of
locally strongly convex hypersurfaces $M$ in $\mathbb{A}^{n+1}$, $n \geq 3$, 
having at every point three distinct affine principal curvatures 
$\lambda _{0}$, $\lambda _{1}$ and $\lambda _{2}$ with multiplicities $1$, $n_{1}$ and $n_{2}$,
respectively. 
Examples of such hypersurfaces are given in \cite{{2014_OB_LV}, {Cece_Li}} (see Example 7.1).
If some additional conditions are satisfied then the tensor $R^{\ast}$ of such hypersurfaces 
is a linear combination of Kulkarni-Nomizu products formed by the tensors $h$, 
$\mathrm{Ric} (R^{\ast})$ and $(\mathrm{Ric} (R^{\ast}))^{2}$ (see Theorem 7.3 for details).
It is clear that locally strongly convex affine $2$-quasi-umbilical hypersurfaces 
belong to this class of hypersurfaces. 

In the last section we apply results from Section 7 to obtain a curvature property of pseudosymmetry type 
of hypersurfaces $M$ isometrically immersed in a Riemannian space forms $N^{n+1}(c)$, $n \geq 3$,
having at every point three distinct principal curvatures $\lambda _{0}$, $\lambda _{1}$ and $\lambda _{2}$, 
with multiplicities $1$, $n_{1}$ and $n_{2}$, respectively. 
Evidently, $2$-quasi-umbilical hypersurfaces belong to this class of hypersurfaces. 
If some additional conditions are satisfied then the curvature tensor $R$ of $M$ 
is a linear combination of Kulkarni-Nomizu products formed by the metric tensor $g$, 
the Ricci tensor $S$ and the tensor $S^{2}$ of $M$ (see Theorem 8.1 for details).
We note that curvature properties of pseudosymmetry type of hypersurfaces 
isometrically immersed in an Euclidean space ${\mathbb{E}}^{n+1}$, $n \geq 5$, 
having at every point three distinct principal curvatures $\lambda _{0}$, $\lambda _{1}$ and $\lambda _{2}$, 
with multiplicities $1$, $n_{1}$ and $n_{2}$, respectively, such that $n_{1} = n_{2} \geq 2$,
were obtained in {\cite[Section 4] {Saw-2015}}. 
In that paper also a class of such hypersurfaces was determined {\cite[Example 5.1] {Saw-2015}}.

\section{Preliminaries}
%\hspace*{6mm}
Throughout this paper, all manifolds are assumed to be connected paracompact
manifolds of class $C^{\infty }$. Let $(M,g)$ be an $n$-dimensional, $n \geq 3$,
semi-Riemannian manifold, and let $\nabla$ be its Levi-Civita connection and $\Xi (M)$ the Lie
algebra of vector fields on $M$. We define on $M$ the endomorphisms 
$X \wedge _{A} Y$ and ${\mathcal{R}}(X,Y)$ of $\Xi (M)$ by
\begin{eqnarray*}
(X \wedge _{A} Y)Z 
&=& 
A(Y,Z)X - A(X,Z)Y, \\ 
{\mathcal R}(X,Y)Z 
&=& 
\nabla _X \nabla _Y Z - \nabla _Y \nabla _X Z - \nabla _{[X,Y]}Z ,
\end{eqnarray*}
respectively, where $A$ is a symmetric $(0,2)$-tensor on $M$ and $X, Y, Z \in \Xi (M)$. 
The Ricci tensor $S$, the Ricci operator ${\mathcal{S}}$ and the scalar curvature
$\kappa $ of $(M,g)$ are defined by 
\begin{eqnarray*}
S(X,Y)\ =\ \mathrm{tr} \{ Z \rightarrow {\mathcal{R}}(Z,X)Y \} ,\ \ \ 
g({\mathcal S}X,Y)\ =\ S(X,Y) ,\ \ \ 
\kappa \ =\ \mathrm{tr}\, {\mathcal{S}},
\end{eqnarray*}
respectively. 
The endomorphism ${\mathcal{C}}(X,Y)$ is defined by
\begin{eqnarray*}
{\mathcal C}(X,Y)Z  &=& {\mathcal R}(X,Y)Z 
- \frac{1}{n-2}(X \wedge _{g} {\mathcal S}Y + {\mathcal S}X \wedge _{g} Y
- \frac{\kappa}{n-1}X \wedge _{g} Y)Z .
\end{eqnarray*}
Now the $(0,4)$-tensor $G$, the Riemann-Christoffel curvature tensor $R$ and
the Weyl conformal curvature tensor $C$ of $(M,g)$ are defined by
%\begin{eqnarray*}
$G(X_1,X_2,X_3,X_4)  = g((X_1 \wedge _{g} X_2)X_3,X_4)$,
$R(X_1,X_2,X_3,X_4)  = g({\mathcal R}(X_1,X_2)X_3,X_4)$ and $C(X_1,X_2,X_3,X_4)  = g({\mathcal C}(X_1,X_2)X_3,X_4)$,
%\end{eqnarray*}
respectively, where $X_1,X_2,\ldots \in \Xi (M)$.  

Let ${\mathcal B}$ be a tensor field sending any $X, Y \in \Xi (M)$ to a skew-symmetric endomorphism 
${\mathcal B}(X,Y)$, 
and let $B$ be
a $(0,4)$-tensor associated with ${\mathcal B}$ by
\begin{eqnarray}
B(X_1,X_2,X_3,X_4) &=& 
g({\mathcal B}(X_1,X_2)X_3,X_4)\, .
\label{DS5}
\end{eqnarray}
It is well-known that the tensor $B$ is said to be a {\sl{generalized curvature tensor}}  if the
following two conditions are fulfilled $B(X_1,X_2,X_3,X_4) = B(X_3,X_4,X_1,X_2)$ and 
\begin{eqnarray*}
B(X_1,X_2,X_3,X_4) + B(X_2,X_3,X_1,X_4) + B(X_3,X_1,X_2,X_4) &=& 0 . 
\end{eqnarray*}

For ${\mathcal B}$ as above, let $B$ be again defined by (\ref{DS5}). 
We extend the endomorphism ${\mathcal B}(X,Y)$
to a derivation ${\mathcal B}(X,Y) \cdot \, $ of the algebra of tensor fields on $M$,
assuming that it commutes with contractions and 
${\mathcal B}(X,Y) \cdot \, f  = 0$ for any smooth function $f$ on $M$. Now for a $(0,k)$-tensor field $T$,
$k \geq 1$, we can define the $(0,k+2)$-tensor $B \cdot T$ by
\begin{eqnarray*}
& & (B \cdot T)(X_1,\ldots ,X_k,X,Y) \ =\ 
({\mathcal B}(X,Y) \cdot T)(X_1,\ldots ,X_k)\\  
&=& - T({\mathcal{B}}(X,Y)X_1,X_2,\ldots ,X_k)
- \cdots - T(X_1,\ldots ,X_{k-1},{\mathcal{B}}(X,Y)X_k)\, .
\end{eqnarray*}
If $A$ is a symmetric $(0,2)$-tensor then we define the
$(0,k+2)$-tensor $Q(A,T)$ by
\begin{eqnarray*}
& & Q(A,T)(X_1, \ldots , X_k, X,Y) \ =\
(X \wedge _{A} Y \cdot T)(X_1,\ldots ,X_k)\\  
&=&- T((X \wedge _A Y)X_1,X_2,\ldots ,X_k) 
- \cdots - T(X_1,\ldots ,X_{k-1},(X \wedge _A Y)X_k)\, .
\end{eqnarray*}
In this manner we obtain the $(0,6)$-tensors $B \cdot B$ and $Q(A,B)$. Substituting 
${\mathcal{B}} = {\mathcal{R}}$ or ${\mathcal{B}} = {\mathcal{C}}$, $T=R$ or $T=C$ or $T=S$, $A=g$ or $A=S$ in the above formulas, 
we get the tensors $R\cdot R$, $R\cdot C$, $C\cdot R$, $R\cdot S$, $Q(g,R)$, $Q(S,R)$, $Q(g,C)$ and $Q(g,S)$.

For a symmetric $(0,2)$-tensor $E$ and a $(0,k)$-tensor $T$, $k \geq 2$, we
define their Kulkarni-Nomizu product $E \wedge T$ by (see, e.g., {\cite[Section 2] {2013_DGHHY}}) 
\begin{eqnarray*}
& &(E \wedge T )(X_{1}, X_{2}, X_{3}, X_{4}; Y_{3}, \ldots , Y_{k})\\
&=&
E(X_{1},X_{4}) T(X_{2},X_{3}, Y_{3}, \ldots , Y_{k})
+ E(X_{2},X_{3}) T(X_{1},X_{4}, Y_{3}, \ldots , Y_{k} )\\
& &
- E(X_{1},X_{3}) T(X_{2},X_{4}, Y_{3}, \ldots , Y_{k})
- E(X_{2},X_{4}) T(X_{1},X_{3}, Y_{3}, \ldots , Y_{k}) .
\end{eqnarray*}
The tensor $E \wedge T$ will be called the {\sl{Kulkarni-Nomizu}} tensor of $E$ and $T$. 

For a symmetric $(0,2)$-tensor $A$ we denote by ${\mathcal{A}}$ the endomorphism
related to $A$ by 
$g({\mathcal{A}}X,Y) = A(X,Y)$.
The tensors
$A^{p}$, $p = 2, 3, \ldots $, are defined by 
$A^{p}(X,Y) = A^{p-1} ({\mathcal{A}}X, Y)$.

It is obvious that the following tensors are generalized curvature tensors: $R$, $C$ and 
$E \wedge F$, where $E$ and $F$ are symmetric $(0,2)$-tensors. 
We have $G = \frac{1}{2}\, g \wedge g$,
\begin{eqnarray}
C &=& R - \frac{1}{n-2}\, g \wedge S + \frac{\kappa }{(n-2) (n-1)} \, G , 
\label{Weyl}
\end{eqnarray}
and (see, e.g., {\cite[Lemma 2.2(i)] {2013_DGHHY}})
\begin{eqnarray}
(a)\ \
Q(E, E \wedge F) \ =\ - \frac{1}{2}\, Q(F, E \wedge E ) , 
\ \ & &\ \
(b)\ \
E \wedge Q(E,F) \ =\ - \frac{1}{2}\, Q(F, E \wedge E ).
\label{DS7}
\end{eqnarray}
Let
$T$ be a $(0,k)$-tensor, $k = 2, 3, \ldots $.  
The tensor $Q(A,T)$ is called the Tachibana tensor of $A$ and $T$, 
or the Tachibana tensor for short (see, e.g., \cite{DGPSS}). 
By an application of (\ref{DS7})(a) we obtain on $M$ the identities
\begin{eqnarray}
\ \ \ \
Q(g, g \wedge S)\ = \ - Q(S,G)
\ \ \ \mbox{and}
\ \ \ Q(S, g \wedge S)\ = \ - \frac{1}{2}\, Q(g, S \wedge S).  
\label{dghz01}
\end{eqnarray}
Using the tensors $g$, $R$ and $S$ we can define the following $(0,6)$-Tachibana tensors: 
$Q(S,R)$, $Q(g,R)$, $Q(g,g \wedge S)$ and $Q(S, g \wedge S)$. 
We can check, by making use of (\ref{DS7})(a) and (\ref{dghz01}), that other $(0,6)$-Tachibana tensors  
constructed from $g$, $R$ and $S$ may be expressed by the four Tachibana tensors mentioned above 
or vanish identically on $M$.
We also have 
\begin{prop}
$\mathrm{ ( }$see, e.g., {\cite[Proposition 2.4] {2013_DGHHY}}$\mathrm{ ) }$  
Let $(M,g)$, $n \geq 3$, be a semi-Rieman\-nian
manifold. Let a non-zero symmetric 
$(0,2)$-tensor $A$ and a generalized curvature tensor $B$, defined at $x \in M$,
satisfy at this point $Q(A,B) = 0$.
In addition, let $Y$ be a vector at $x$ such that the scalar $\rho = w(Y)$ is non-zero,
where $w$ is a covector defined by $w(X) = A(X,Y)$, $X \in T_{x}M$. 
Then we have:
\newline 
(i) $A - \rho \, w \otimes w \neq 0$ and $B = \lambda \, A \wedge A$, $\lambda \in {\mathbb R}$, 
\newline
(ii) $A = \rho \,  w \otimes w$ and 
\begin{eqnarray*}
w(X)\, B(Y,Z, X_{1}, X_{2}) + w(Y)\, B(Z,X, X_{1}, X_{2}) + w(Z)\, B(X,Y, X_{1}, X_{2}) &=& 0, 
%\label{TT00}
\end{eqnarray*}
where $X,Y,Z,X_{1},X_{2} \in T_{x}M$.
\newline
Moreover, in both cases the following condition holds at $x$
\begin{eqnarray*}
B \cdot B &=& Q( \mathrm{Ric}(B),B). 
%\label{TT01}
\end{eqnarray*}
\end{prop}

\section{Some special generalized curvature tensors}

Let  
$\{ e_{1}, e_{2}, \ldots , e_{n} \}$ be an orthonormal basis 
of $T_{x}M$  at a point $x \in M$ of a semi-Riemannian manifold
$(M,g)$, $n \geq 3$, and let
$g( e_{j}, e_{k}) = \varepsilon _{j} \delta _{jk}$, where
$\varepsilon _{j} = \pm 1$ and $h,i,j,k,l,m, r,s \in \{ 1, 2, \ldots , n \}$.
For a generalized curvature tensor $B$ on $M$ we denote by 
$\mathrm{Ric}(B)$, $\kappa (B)$ and $\mathrm{Weyl}(B)$   
its Ricci tensor, the scalar curvature and the Weyl tensor, 
respectively. 
We have
\begin{eqnarray}
\mathrm{Ric}(B)(X,Y)
&=&
\sum _{j=1}^{n}
\varepsilon_{j}\, B( e_{j}, X,Y, e_{j}),\ \ \
\kappa(B) \ =\ \sum _{j=1}^{n}
\varepsilon_{j}\, \mathrm{Ric}(B)( e_{j}, e_{j}),\nonumber\\
Weyl(B)
&=&
B - \frac{1}{n-2}\, g \wedge \mathrm{Ric}(B) +
\frac{\kappa(B)}{ (n-2)(n-1)}\, G .
\label{AAconh01}
\end{eqnarray}
We define the following subsets of $M$:
${\mathcal{U}}_{ B} =
\{ x \in M\, | \, B \neq  ( \kappa (B) / ( (n-1) n) )\, G \ \mbox{at}\ x \}$,
${\mathcal{U}}_{ \mathrm{Ric}(B)} =
\{ x \in M\, | \, \mathrm{Ric}(B) \neq ( \kappa (B) / n )\, g \ \mbox{at}\ x \}$ and
${\mathcal{U}}_{ \mathrm{Weyl}(B)} =
\{ x \in M\, | \, \mathrm{Weyl}(B) \neq 0\ \mbox{at}\ x \}$.
We note that
${\mathcal{U}}_{B} = 
{\mathcal{U}}_{ \mathrm{Ric}(B)} \cup {\mathcal{U}}_{ \mathrm{Weyl}(B)}$  (cf. \cite{2013_DGHHY}).  

Let $B_{hijk}$, $T_{hijk}$ and $A_{ij}$ be the local components of the generalized curvature tensors $B$ and $T$ 
and a symmetric $(0,2)$-tensor $A$ on $M$, respectively. 
The local components $(B \cdot T)_{hijklm}$ and $Q(A,T)_{hijklm}$ of the tensors $B \cdot T$ and $Q(A,T)$ 
are the following:
\begin{eqnarray}
(B \cdot T)_{hijklm} 
&=&
g^{rs}(
T_{rijk}B_{ shlm}
+ T_{hrjk}B_{ silm}
+ T_{hirk}B_{sjlm}
+ T_{hijr}B_{sklm}) ,
\label{D01}\\
%\end{eqnarray*}
%\begin{eqnarray*}
Q(A,T)_{hijklm}
&=&
A_{hl}T_{ mijk} 
+ A_{il}T_{ hmjk} 
+ A_{jl}T_{ himk} 
+ A_{kl}T_{ hijm}\nonumber\\ 
& &-
A_{hm}T_{ lijk} 
- A_{im}T_{ hljk} 
- A_{jm}T_{ hilk} 
- A_{km}T_{ hijl}  .
\label{D02}
\end{eqnarray}

Let $B$ be a generalized curvature tensor on a semi-Rieman\-nian manifold
$(M,g)$, $n \geq 4$. 
The local components $(B \cdot \mathrm{Weyl}(B))_{hijklm}$ and $( \mathrm{Weyl}(B) \cdot B)_{hijklm}$ of the tensors 
$B \cdot \mathrm{Weyl}(B)$ and $\mathrm{Weyl}(B) \cdot B$ 
are the following:
\begin{eqnarray}
(B \cdot \mathrm{Weyl}(B))_{hijklm} 
&=&
g^{rs}(
\mathrm{Weyl}(B)_{rijk}B_{ shlm}
+ \mathrm{Weyl}(B)_{hrjk}B_{ silm}\nonumber\\
& &
+ \mathrm{Weyl}(B)_{hirk}B_{sjlm}
+ \mathrm{Weyl}(B)_{hijr}B_{sklm}) ,
\label{gencurv75}
\end{eqnarray}
\begin{eqnarray}
(\mathrm{Weyl}(B) \cdot B)_{hijklm} 
&=&
g^{rs}(
B_{rijk} \mathrm{Weyl}(B)_{ shlm}
+ B_{hrjk} \mathrm{Weyl}(B)_{ silm}\nonumber\\
& &
+ B_{hirk} \mathrm{Weyl}(B)_{sjlm}
+ B_{hijr} \mathrm{Weyl}(B)_{sklm}) .
\label{gencurv74}
\end{eqnarray}
We set $V_{mijk} = g^{rs} Ric(B)_{mr}B_{sijk}$. 
Now (\ref{AAconh01}), (\ref{gencurv75}) and (\ref{gencurv74}) give 
(cf. \cite{{DGHSaw}})
\begin{eqnarray*}
& &
(n-2)(B \cdot \mathrm{Weyl}(B)  - \mathrm{Weyl}(B) \cdot B)_{hijklm} \ =\ Q( \mathrm{Ric}(B),B)_{hijklm}\nonumber\\
& & 
- \frac{\kappa}{n-1}\, Q(g,B)_{hijklm} + g_{hl}V_{mijk} - g_{hm}V_{lijk} - g_{il}V_{mhjk} + g_{im}V_{lhjk}\nonumber\\
& &
+ g_{jl}V_{mkhi} - g_{jm}V_{lkhi} - g_{kl}V_{mjhi} + g_{km}V_{ljhi} 
- g_{ij}(B\cdot \mathrm{Ric}(B))_{hklm}\nonumber\\
& &
- g_{hk}(B \cdot \mathrm{Ric}(B))_{ijlm} + g_{ik}(B \cdot \mathrm{Ric}(B))_{hjlm} + g_{hj} (B \cdot \mathrm{Ric}(B))_{iklm} .
%\label{gencurv72}
\end{eqnarray*}

\begin{prop} {\cite[Lemma 3.3, Theorem 3.1] {Kow2}} 
Let $B$ be a $(0,4)$-tensor on a semi-Riemannian manifold $(M,g)$, $n \geq 3$,
defined by
\begin{eqnarray}
B &=& \frac{\phi }{2}\, A \wedge A + \mu \, g \wedge A + \frac{\eta}{2}\, g \wedge g ,
\label{aff0303}
\end{eqnarray}
where $A$ is a symmetric $(0,2)$-tensor 
and $\phi, \mu, \eta$ are some functions defined on $M$. 
\newline
(i) 
Let $\mathcal{U}$ be the set of all points of $M$ at which $\phi$ is non-zero.
Then the following conditions are satisfied on $\mathcal{U}$
\begin{eqnarray*}
A^{2} &=& 
\phi ^{-1} ( 
(\phi \mathrm{tr}(A) + (n-2) \mu )\, A  
+ (\mu \mathrm{tr}(A) + (n-1) \eta )\, g 
- \mathrm{Ric} (B) ),
%\label{aff21}
\\
\ \ \ \ \ \ \
B \cdot A &=& 
Q(  \mathrm{Ric}(B) + (n-2) (\mu ^{2} - \phi \eta ) \phi^{-1} \, g, A + \mu \phi ^{-1}\, g ) ,
%\label{aff22} 
\\
B \cdot B &=& Q(\mathrm{Ric}(B), B) + (n-2) (\mu ^{2} - \phi \eta ) \phi^{-1} \, Q(g, \mathrm{Weyl} (B)) . 
%\label{aff23}
\end{eqnarray*}
(ii) If $\phi = 0$ at a point of $M$ then $\mathrm{Weyl} (B)$ vanishes at this point. 
\end{prop}

%\vspace{3mm}

According to \cite{P106}, 
a generalized curvature tensor $B$ on a semi-Riemannian ma\-ni\-fold $(M,g)$, $n \geq 4$,
is called a {\sl{Roter type tensor}} if 
\begin{eqnarray}
B &=& \frac{\phi}{2}\, \mathrm{Ric}(B) \wedge \mathrm{Ric}(B) + \mu\, g \wedge \mathrm{Ric}(B) + \frac{\eta }{2}\, g \wedge g 
\label{eq:h7}
\end{eqnarray}
on ${\mathcal{U}}_{\mathrm{Ric}(B)} \cap {\mathcal{U}}_{\mathrm{Weyl} (B)}$,
whereby $\phi$, $\mu $ and $\eta $ are some functions on this set.
Evidently, (\ref{eq:h7}) is a special case of (\ref{aff0303}).
Manifolds admitting 
Roter type tensors were investigated among others in \cite{{2013_DGHHY}, {2015_DGHZ}, {Kow2}}. 
We have 
\begin{prop}
Let $B$ be a generalized curvature tensor on a semi-Rieman\-nian manifold $(M,g)$, $n \geq 4$,
satisfying (\ref{eq:h7}) on 
${\mathcal{U}} = {\mathcal{U}}_{Ric(B)} \cap {\mathcal{U}}_{Weyl (B)} \subset M$.
\newline
(i) 
$\mathrm{ ( }$cf. {\cite[Section 3] {2015_DGHZ}}, {\cite[Sections 1 and 4] {Kow2}}$\mathrm{ ) }$
The following relations hold on ${\mathcal{U}}$ 
\begin{eqnarray}
&(a)&\ \ \ (\mathrm{Ric}(B))^{2} 
\ =\ \alpha _{1}\, \mathrm{Ric}(B) + \alpha _{2} \, g ,\nonumber\\
&(b)&\ \ \
\alpha _{1} \ =\ \kappa(B) + \phi^{-1} ( (n-2)\mu -1 ),\ \ \ \alpha _{2}
\ =\
\phi^{-1} (  \mu \kappa(B) + (n-1) \eta ) ,
\label{roter71}
\end{eqnarray}
\begin{eqnarray}
&(a)&\ \ \ B \cdot B \ =\ L_{B}\, Q(g,B),
\ \ \
L_{B}
\ =\
\phi^{-1}  \left( (n-2) (\mu ^{2} - \phi \eta) - \mu \right),\nonumber\\
&(b)&\ \ \ B \cdot \mathrm{Weyl}(B) \ =\ L_{B}\, Q(g,\mathrm{Weyl}(B)),\nonumber\\
&(c)&\ \ \ B \cdot B \ =\ Q(\mathrm{Ric}(B),B) + 
L\, Q(g,\mathrm{Weyl}(B)) ,\ \ \
L \ = \ L_{B} + \phi^{-1}  \mu , 
\label{roter72}
\end{eqnarray}
\begin{eqnarray}
\mathrm{Weyl}(B) \cdot B &=&  L_{\mathrm{Weyl}(B)}\, Q(g,B) ,\ \ \ 
L_{\mathrm{Weyl}(B)} \ =\  L_{B} + \frac{1}{n-2} \left( \frac{\kappa(B) }{n-1} - \alpha _{1} \right) ,
\label{roter73NN}
\end{eqnarray}
\begin{eqnarray*}
\mathrm{Weyl}(B) \cdot \mathrm{Weyl}(B) &=& 
L_{ \mathrm{Weyl}(B)}\, 
Q(g,\mathrm{Weyl}(B)) ,
%\label{roter73}
\end{eqnarray*}
\begin{eqnarray}
B \cdot \mathrm{Weyl}(B) - \mathrm{Weyl}(B) \cdot B &=& 
\left( \phi^{-1}  ( \mu - \frac{1}{n-2} ) + \frac{\kappa(B) }{n-1} \right) Q(g,B)\nonumber\\
& &+ \left(  \phi^{-1} \mu ( \mu - \frac{1}{n-2}) - \eta \right) Q(\mathrm{Ric}(B),G) .
\label{dh3}
\end{eqnarray}
(ii) {\cite[Proposition 3.2 (ii)] {2015_DGHZ}}
The following relation holds on ${\mathcal{U}}$ 
\begin{eqnarray}
Q(\mathrm{Ric}(B),\mathrm{Weyl}(B)) &=&   
\phi ^{- 1} ( \frac{1}{n-2} - \mu)\, Q(g,B)\nonumber\\
& &
+ 
\frac{1}{n-2} \left( L_{B} - \frac{\kappa(B) }{n-1} \right)  Q(g, g \wedge \mathrm{Ric}(B)) .
\label{roter74}
\end{eqnarray}
(iii) $\mathrm{ ( }$cf. {\cite[Section 3] {2016_DGHZhyper}}, {\cite[Proposition 3.3] {2016_DGJZ}}$\mathrm{ ) }$
The following relation holds on ${\mathcal{U}}$ 
\begin{eqnarray} 
\mathrm{Weyl}(B) \cdot B - B \cdot \mathrm{Weyl}(B) 
&=& Q(\mathrm{Ric}(B),\mathrm{Weyl}(B)) - \frac{\kappa(B)}{n-1}\, Q( g , \mathrm{Weyl}(B)) . 
\label{aff024}  
\end{eqnarray}
(iv)
The following relation holds on ${\mathcal{U}}$ 
\begin{eqnarray} 
\mathrm{Weyl}(B) \cdot B + B \cdot \mathrm{Weyl}(B) 
&=& Q(\mathrm{Ric}(B),\mathrm{Weyl}(B))\nonumber\\
& & 
+
\left( L + L_{ \mathrm{Weyl}(B)} - \frac{1}{(n-2) \phi } \right) Q( g , \mathrm{Weyl}(B)) . 
\label{aff024777}  
\end{eqnarray}
(v)
The following relations hold on ${\mathcal{U}}$ 
\begin{eqnarray}
\alpha _{2} &=& \frac{1}{n}\, ( \mathrm{tr}_{g} ( (\mathrm{Ric} (B))^{2} )  - \alpha _{1} \kappa (B) ) ,
\label{aff024anew} \\
(\mathrm{Ric}(B))^{2} - \frac{1}{n}\, \mathrm{tr}_{g} ((\mathrm{Ric} (B))^{2})\, g 
&=& 
\alpha _{1}\, ( \mathrm{Ric}(B) - \frac{1}{n} \kappa (B) \, g ) ,
\label{aff024bnew} \\
(\mathrm{Ric}(B))^{3} - \frac{1}{n}\, \mathrm{tr}_{g} ( (\mathrm{Ric} (B))^{2} )\, \mathrm{Ric}(B) 
&=& 
\alpha _{1}\,  ( (\mathrm{Ric}(B))^{2} - \frac{1}{n} \kappa (B) \, \mathrm{Ric}(B) ) .
\label{aff024cnew} 
\end{eqnarray}
\end{prop}
{\bf{Proof.}} (iii)
The condition (\ref{aff024}) 
is an immediate consequence of 
(\ref{AAconh01}), (\ref{roter72})(a), (\ref{dh3}) and (\ref{roter74}).
(iv)
The condition (\ref{aff024}), 
by (\ref{roter72})(b),
turns into
\begin{eqnarray*} 
\mathrm{Weyl}(B) \cdot B + B \cdot \mathrm{Weyl}(B) 
&=& Q(\mathrm{Ric}(B),\mathrm{Weyl}(B)) 
+ \left( 2 L_{B} - \frac{\kappa(B)}{n-1} \right) Q( g , \mathrm{Weyl}(B)) . 
\end{eqnarray*}
Using now
(\ref{roter71})(b),
(\ref{roter72})(c) and
(\ref{roter73NN})
we can easily check that
\begin{eqnarray*} 
L + L_{ \mathrm{Weyl}(B)} - \frac{1}{(n-2) \phi } &=& 2 L_{B} - \frac{\kappa(B)}{n-1}.
\end{eqnarray*}
The last remark completes the proof of (iv).
(v) 
Contracting (\ref{roter71})(a) we get (\ref{aff024anew}). 
Now (\ref{roter71})(a), by (\ref{aff024anew}), turns into 
(\ref{aff024bnew}). 
From (\ref{aff024bnew}), by a suitable contraction, we get easily (\ref{aff024cnew}),
completing the proof.
\qed

\vspace{2mm}

Evidently,  (\ref{aff08}) is a special case of (\ref{aff024}).  We also have 
\begin{prop}
Let $B$ be a generalized curvature tensor on a semi-Rieman\-nian manifold $(M,g)$, $n \geq 4$.
Let ${\mathcal{U}} \subset {\mathcal{U}}_{Ric(B)} \cap {\mathcal{U}}_{Weyl (B)} \subset M$
be the set of all points at which 
the tensor $( \mathrm{Ric}(B) )^{2}$ 
is not a linear combination of the tensors
$g$ and $\mathrm{Ric}(B)$. 
\newline
(i) If the condition 
\begin{eqnarray}
\ \ \ \ 
\mathrm{Weyl}(B) &=& \frac{\phi }{2}\, \mathrm{Ric}(B) \wedge \mathrm{Ric}(B)
+ \alpha _{1}\,  g \wedge ( \mathrm{Ric}(B) )^{2}
+ \alpha _{2} \,  g \wedge \mathrm{Ric}(B) 
+ \frac{\alpha _{3} }{2}\, g \wedge g 
\label{aff02424aff}
\end{eqnarray}
is satisfied on ${\mathcal{U}}$, 
where $\alpha _{1} , \alpha _{2} , \alpha _{3}$ are some functions on this set,
then 
\begin{eqnarray}
\alpha _{1} \ =\ \frac{\phi }{n-2}\, ,\ \ \alpha _{2} \ =\ - \frac{\kappa(B) \phi }{n-2}\, ,\ \ 
\alpha _{3} \ =\ \frac{ ( (\kappa(B))^{2} - \mathrm{tr}( ( \mathrm{Ric}(B) )^{2})     ) \phi }{(n-2)(n-1)} .
\label{aff03535aff}
\end{eqnarray}
(ii) 
If the condition 
\begin{eqnarray}
\ \ \ \ 
B &=& \frac{\phi }{2}\, \mathrm{Ric}(B) \wedge \mathrm{Ric}(B)
+ \beta _{1}\,  g \wedge ( \mathrm{Ric}(B) )^{2}
+ \beta  _{2} \,  g \wedge \mathrm{Ric}(B) 
+ \frac{\beta _{3} }{2}\, g \wedge g 
\label{raff02424raff}
\end{eqnarray}
is satisfied on ${\mathcal{U}}$, 
where $\beta _{1} , \beta _{2} , \beta _{3}$ are some functions on this set,
then 
\begin{eqnarray}
\ \ \ \ \ \ 
\beta _{1} \ =\ \frac{\phi }{n-2}\, ,\ \ \beta _{2} \ =\ \frac{1 - \kappa(B) \phi }{n-2}\, ,\ \ 
\beta _{3} \ =\ \frac{  ( (\kappa(B))^{2} - \mathrm{tr}( ( \mathrm{Ric}(B) )^{2})     ) \phi  - \kappa (B) }{(n-2)(n-1)} .
\label{raff03535raff}
\end{eqnarray}
\end{prop}
{\bf{Proof.}} 
(i) The condition (\ref{aff02424aff}), by a suitable contraction
and making use of our assumptions, yields (\ref{aff03535aff}).
(ii) The condition (\ref{raff02424raff}), by (\ref{AAconh01}),
turns into
\begin{eqnarray*}
\ \ \ \ 
Weyl (B) &=& \frac{\phi }{2}\, \mathrm{Ric}(B) \wedge \mathrm{Ric}(B)
+ \beta _{1}\,  g \wedge ( \mathrm{Ric}(B) )^{2}
+  (\beta  _{2} - \frac{1}{n-2}) \,  g \wedge \mathrm{Ric}(B)\\
& & 
+ \left( \beta _{3} + \frac{ \kappa (B)}{(n-2)(n-1)} \right) \frac{1}{2}\, g \wedge g .
%\label{raff02424raff}
\end{eqnarray*}
From this, in view of (i), we get immediately (\ref{raff03535raff}).
\qed

\begin{prop} {\cite[Lemma 2.1] {2016_DGJZ}} 
Let $(M,g)$, $n \geq 3$, be a semi-Riemannian manifold.
Let $A$ be a symmetric $(0,2)$-tensor on $M$ such that 
$\mathrm{rank}(A) = 2$ at some point $x \in M$.
\newline
(i) cf. {\cite[Lemma 2.1] {P106}} The tensors $A$, $A^{2}$ and $A^{3}$ satisfy at $x$ the following relations
\begin{eqnarray}
A^{3} &=& \mathrm{tr}(A)\, A^{2} + \frac{1}{2} ( \mathrm{tr}(A^{2}) - (\mathrm{tr}(A))^{2})\, A ,
\label{eqn10.1}\\ 
A \wedge A^{2} &=& \frac{1}{2} \mathrm{tr}(A)\, A \wedge A ,
\label{eqn14.1}\\
A^{2} \wedge A^{2} &=& - \frac{1}{2} ( \mathrm{tr}(A^{2}) - (\mathrm{tr}(A))^{2})\, A \wedge A ,
\label{eqn14.2}
\end{eqnarray}
\begin{eqnarray}
(A^{2} - \mathrm{tr}(A)\, A) \wedge  (A^{2} - \mathrm{tr}(A)\, A )
&=& - \frac{1}{2} ( \mathrm{tr}(A^{2}) - (\mathrm{tr}(A))^{2})\, A \wedge A .
\label{eqn14.2dd}
\end{eqnarray}
(ii) Let $T$ be a generalized curvature tensor on $M$ satisfying 
\begin{eqnarray*}
T &=& \frac{\phi _{0}}{2}\, A \wedge A + \phi _{2}\, g \wedge A + \phi _{3}\, G + \phi _{4}\, g \wedge A^{2}
+ \phi _{5}\, A \wedge A^{2} + \frac{ \phi _{6}}{2}\, A^{2} \wedge A^{2} , 
%\label{eqn14.3}
\end{eqnarray*}
where $\phi _{0}$, $\phi _{2},\ \ldots \ , \phi _{6}$
are some functions on $M$. Then at given point $x$ we have  
\begin{eqnarray*}
T &=& \frac{\phi _{1}}{2}\, A \wedge A + \phi _{2}\, g \wedge A + \phi _{3}\, G + \phi _{4}\, g \wedge A^{2} ,\\
\phi _{1} &=& \phi _{0} +  \mathrm{tr}(A)\, \phi _{5}
- \frac{1}{2} ( \mathrm{tr}(A^{2}) - (\mathrm{tr}(A))^{2})\, \phi _{6} .
\end{eqnarray*}
\end{prop}

\section{Pseudosymmetry type curvature conditions}

It is well-known that if a semi-Riemannian manifold $(M,g)$, $n \geq 3$, 
is locally symmetric then $\nabla R = 0$ on $M$. This  
implies  the following integrability condition
${\mathcal{R}}(X,Y ) \cdot R = 0$, in short 
\begin{eqnarray}
R \cdot R &=& 0 .
\label{semisymmetry}
\end{eqnarray}
Semi-Riemannian manifolds satisfying (\ref{semisymmetry})
are called {\sl semisymmetric}, see, e.g., 
{\cite[Chapter 8.5.3] {M7}} and 
{\cite[Chapter 1.6] {M6}}.
Semisymmetric manifolds form a subclass of the class of pseudosymmetric manifolds.
A semi-Riemannian manifold $(M,g)$, $n \geq 3$, is said to be {\sl pseudosymmetric} 
if the tensors $R \cdot R$ and $Q(g,R)$ are linearly dependent at every point of $M$
(see, e.g., {\cite[Chapter 8.5.3] {M7}}, {\cite[Chapter 12.4] {M6}}, 
\cite{{2013_DGHHY}, {DGHSaw}, {DHV2008}, {DVV1991}, {SDHJK}} 
and references therein). 
This is equivalent to
\begin{eqnarray}
R \cdot R &=& L_{R}\, Q(g,R) 
\label{pseudo}
\end{eqnarray}
on  the set ${\mathcal{U}}_{R} = \{x \in M\, | \, R - \frac{\kappa }{(n-1)n}\, G \neq 0\ \mbox {at}\ x \}$,
where $L_{R}$ is some function on this set. 
Every semisymmetric manifold is pseudosymmetric.
The converse statement is not true (see, e.g., \cite{{DVV1991}}).
We note that (\ref{pseudo}) implies
\begin{eqnarray}
(a)\ \ R \cdot S \ =\ L_{R}\, Q(g,S)\ \ &\mbox{and}& \ \ (b)\ \ R \cdot C \ =\  L_{R}\, Q(g,C) .
\label{Weyl-pseudo-bis}
\end{eqnarray}

It is well-known that $(M,g)$ is said to be an {\sl Einstein manifold} if at every point of $M$ 
its Ricci tensor $S$ is proportional to the metric tensor $g$, 
i.e., $S = \frac{\kappa}{n}\, g$ on $M$. 
We denote by ${\mathcal U}_{S}$ the set of all points of $(M,g)$ at which 
$S$ is not proportional to $g$, i.e.,
${\mathcal U}_{S} \, = \,  \{x \in M\, | \, 
S - \frac{\kappa }{n}\, g \neq 0\ \mbox {at}\ x \}$ and by
${\mathcal U}_{C}$ the set of all points of $M$ at which $C \neq 0$.
We note that
${\mathcal{U}}_{S} \cup {\mathcal{U}}_{C} = {\mathcal{U}}_{R}$ (see, e.g., \cite{2013_DGHHY}). 
The conditions 
(\ref{pseudo}), (\ref{Weyl-pseudo-bis})(a) and (\ref{Weyl-pseudo-bis})(b)
are equivalent on the set ${\mathcal{U}}_{S} \cap {\mathcal{U}}_{C}$ 
of any warped product manifold $M_{1} \times_{F} M_{2}$, 
with $\dim\, M_{1} = \dim M_{2} = 2$ 
(see, e.g., \cite{{2016_DGJZ}} and references therein).
We mention that the Schwarzschild spacetime, the Kottler spacetime, the Reissner-Nordstr\"{o}m spacetime, 
as well as the Friedmann-Lema{\^{\i}}tre-Robertson-Walker spacetimes are the "oldest" examples 
of pseudosymmetric warped product manifolds (see, e.g., \cite{{2016_DGJZ}, {DHV2008}, {DVV1991}, {SDHJK}}).

A semi-Riemannian manifold $(M,g)$, $n \geq 3$, is called {\sl Ricci-pseudosymmetric} 
if the tensors $R \cdot S$ and $Q(g,S)$ are linearly dependent at every point of $M$
(see, e.g., {\cite[Chapter 8.5.3] {M7}}, \cite{DGHSaw}).
This is equivalent on ${\mathcal{U}}_{S}$ to 
\begin{eqnarray}
R \cdot S &=& L_{S}\, Q(g,S) , 
\label{Riccipseudo07}
\end{eqnarray}
where $L_{S}$ is some function on this set. 
Every warped product manifold $\overline{M} \times _{F} \widetilde{N}$
with a $1$-dimensional $(\overline{M}, \overline{g})$ manifold and
an $(n-1)$-dimensional Einstein semi-Riemannian manifold $(\widetilde{N}, \widetilde{g})$, $n \geq 3$, 
and a warping function $F$, 
is a Ricci-pseudosymmetric manifold,
see, e.g., {\cite[Section 1] {Ch-DDGP}} and {\cite[Example 4.1] {2016_DGJZ}}.

A semi-Riemannian manifold $(M,g)$, $n \geq 4$, is said to be {\sl Weyl-pseudo\-sym\-met\-ric} 
if the tensors $R \cdot C$ and $Q(g,C)$ are linearly dependent at every point of $M$
\cite{{2013_DGHHY}, {DGHSaw}}. 
This is equivalent on ${\mathcal{U}}_{C}$ to 
\begin{eqnarray}
R \cdot C &=& L_{1}\, Q(g,C) ,  
\label{Weyl-pseudo}
\end{eqnarray}
where $L_{1}$ is some function on this set. 
Using (\ref{Weyl}), we can check that 
on every Einstein manifold $(M,g)$, $n \geq 4$,
(\ref{Weyl-pseudo}) turns into
\begin{eqnarray*}
R \cdot R &=& L_{1}\, Q(g,R) .
\end{eqnarray*}
For a presentation of results on the problem 
of the equivalence of pseudosymmetry, Ricci-pseudosymmetry and Weyl-pseudosymmetry
we refer to {\cite[Section 4] {DGHSaw}}.  

A semi-Riemannian manifold $(M,g)$, $n \geq 4$, is said to have {\sl pseudosymmetric Weyl tensor}
if the tensors $C \cdot C$ and $Q(g,C)$ are linearly dependent at every point of $M$ 
(see, e.g., \cite{{2013_DGHHY}, {DGHSaw}, {2016_DGJZ}}).
This is equivalent on ${\mathcal U}_{C}$ to 
\begin{eqnarray}
C \cdot C &=& L_{C}\, Q(g,C) ,  
\label{4.3.012}
\end{eqnarray}
where $L_{C}$ is some function on this set. 
Every warped product manifold $M_{1} \times_{F} M_{2}$, 
with $\dim\, M_{1} = \dim M_{2} = 2$, sa\-tis\-fies (\ref{4.3.012})
(see, e.g., \cite{{2013_DGHHY}, {DGHSaw}, {2016_DGJZ}} and references therein).
Thus in particular,
the Schwarzschild spacetime, the Kottler spacetime
and the Reissner-Nordstr\"{o}m spacetime satisfy (\ref{4.3.012}).
Recently semi-Riemannian manifolds with pseudosymmetric Weyl tensor 
were investigated in \cite{{2013_DGHHY}, {DeHoJJKunSh}}.

Warped product manifolds $\overline{M} \times _{F} \widetilde{N}$, of dimension $\geq 4$,
satisfying on 
${\mathcal U}_{C} \subset \overline{M} \times \widetilde{N}$,
the condition 
\begin{eqnarray}
R \cdot R - Q(S,R) &=& L\, Q(g,C) ,  
\label{genpseudo01}
\end{eqnarray}
where $L$ is some function on this set,
were studied among others in \cite{49}.
In \cite{49} necessary and sufficient conditions for  
$\overline{M} \times _{F} \widetilde{N}$ to be a manifold satisfying (\ref{genpseudo01}) are given.
Moreover, in that paper it was proved that 
any $4$-dimensional warped product manifold $\overline{M} \times _{F} \widetilde{N}$, 
with a $1$-dimensional base $(\overline{M},\overline{g})$, 
satisfies (\ref{genpseudo01}) {\cite[Theorem 4.1] {49}}.

We refer to
\cite{{Ch-DDGP}, {1995_DDDVY}, {49}, {2013_DGHHY}, {DGHSaw}, {2015_DGHZ}, 
{2016_DGJZ}, {2018_DGZ}, {DHV2008}, {DeHoJJKunSh}, {Saw-2015}, {SDHJK}} 
for details on semi-Riemannian manifolds satisfying 
(\ref{pseudo}) and (\ref{Riccipseudo07})-(\ref{genpseudo01}), 
as well other conditions of this kind, named pseudosymmetry type curvature conditions. 
We also refer to {\cite[Section 3] {DeHoJJKunSh}} for a recent survey on manifolds 
satisfying such curvature conditions.
It seems that the condition (\ref{pseudo}) 
is the most important condition of that family of curvature conditions
(see, e.g., \cite{2016_DGJZ}).

%\vspace{3mm}

\section{Roter spaces}

The manifold $(M,g)$ is said to be 
a {\sl quasi-Einstein manifold} if 
\begin{eqnarray}
\mathrm{rank}(S - \alpha\, g) &=& 1
\label{quasi02}
\end{eqnarray}
on ${\mathcal U}_{S} \subset M$, where $\alpha $ is some function on ${\mathcal U}_{S}$.
It is known that every warped product manifold $\overline{M} \times _{F} \widetilde{N}$
with an $1$-dimensional $(\overline{M}, \overline{g})$ base manifold and
a $2$-dimensional manifold $(\widetilde{N}, \widetilde{g})$
or an $(n-1)$-dimensional Einstein manifold
$(\widetilde{N}, \widetilde{g})$, $n \geq 4$, and a warping function $F$,
is a  quasi-Einstein manifold (see, e.g., {\cite{{Ch-DDGP}, {2016_DGJZ}}}). 
We note that (\ref{quasi02}) implies
{\cite[eq. (8)] {G6}}
\begin{eqnarray}
S^{2} &=& ( \kappa - (n-2) \alpha )\, S - \alpha ( \kappa - (n-1) \alpha )\, g  .
\label{quasi0202}
\end{eqnarray}
We mention that quasi-Einstein manifolds arose during the study of exact solutions
of the Einstein field equations and the investigation on quasi-umbilical hypersurfaces 
of conformally flat spaces, see, e.g., \cite{{DGHSaw}, {2016_DGJZ}} and references therein. 
Quasi-Einstein hypersurfaces in semi-Riemannian spaces of constant curvature
were studied among others in
\cite{{R102}, {G6}}, see also \cite{DGHSaw} and references therein.
Quasi-Einstein manifolds satisfying some pseudosymmetry type curvature conditions  
were investigated recently in
\cite{{Ch-DDGP}, {2013_DGHHY}, {2015_DGHZ}, {DeHoJJKunSh}}. 

According to {\cite[Foreword] {Chen-2017}} 
(precisely, {\cite[p. 20] {V2}})  
a Riemannian manifold $(M,g)$, $n \geq 3$, is said to be {\sl{partially Einstein space}} 
if at every point of $M$ its Ricci operator  ${\mathcal{S}}$ has at most two distinct eigenvalues. 
More generally, a semi-Riemannian manifold $(M,g)$, $n \geq 3$, will be called {\sl{partially Einstein space}} 
if at every point $x \in {\mathcal{U}}_{S} \subset M$ its Ricci operator ${\mathcal{S}}$ satisfies
\begin{eqnarray*}
{\mathcal{S}}^{2} &=& \lambda \, {\mathcal{S}} + \mu \, I\!d_{x} , 
\end{eqnarray*}
or equivalently, 
\begin{eqnarray}
S^{2} &=&  \lambda \, S + \mu \, g ,
\label{partiallyEinstein}
\end{eqnarray}
where $\lambda , \mu \in {\mathbb{R}}$  
and $I\!d_{x}$ is the identity transformation of $T_{x} M$.
We can also present (\ref{partiallyEinstein}) in the form (cf. the proof of Proposition 3.2 (iv)) 
\begin{eqnarray*} 
S^{2} - \frac{ \mathrm{tr}_{g} (S^{2})}{n} \, g &=&  \lambda \, ( S - \frac{\kappa }{n} \, g ) . 
\end{eqnarray*}
Evidently, (\ref{quasi0202}) is a special case of (\ref{partiallyEinstein}). 
Thus, every quasi-Einstein manifold is a partially Einstein manifold.
The converse statement is not true.

%\vspace{3mm}

A semi-Riemannian manifold $(M,g)$, $n \geq 3$, 
is called a {\sl $2$-quasi-Einstein manifold} if 
\begin{eqnarray}
\mathrm{rank}(S - \alpha \, g ) &\leq & 2 
\label{quasi0202weak}
\end{eqnarray}
on ${\mathcal U}_{S}$ and $\mathrm{rank} (S - \alpha \, g ) = 2$
on some open non-empty subset of ${\mathcal U}_{S}$, 
where $\alpha $ is some function on ${\mathcal U}_{S}$ 
(see, e.g., \cite{{DGP-TV02}}).
Every warped product manifold $\overline{M} \times _{F} \widetilde{N}$
with a $2$-dimensional base manifold $(\overline{M}, \overline{g})$ 
and a $2$-dimensional manifold $(\widetilde{N}, \widetilde{g})$
or an $(n-2)$-dimensional Einstein semi-Riemannian manifold
$(\widetilde{N}, \widetilde{g})$, when $n \geq 5$, 
and a warping function $F$ satisfies (\ref{quasi0202weak})
(see, e.g., {\cite[Theorem 6.1] {2016_DGJZ}}). 
Thus some exact solutions of the Einstein field equations are 
non-conformally flat $2$-quasi-Einstein manifolds.
For instance, the Reissner-Nordstr\"{o}m spacetime, as well as
the Reissner-Nordstr\"{o}m-de Sitter type spacetimes are such manifolds (see, e.g.,
\cite{Kow2}). 
It seems that the Reissner-Nordstr\"{o}m spacetime
is the "oldest" example of 
a non-conformally flat 
$2$-quasi-Einstein warped product manifold
{\cite[Section 1] {2016_DGJZ}}.
It is known that every $2$-quasi-umbilical hypersurface
in a semi-Riemannian space of constant curvature is a $2$-quasi-Einstein manifold
(see, e.g., \cite{DGP-TV02}).

%\vspace{3mm}

If a non-quasi-Einstein semi-Riemannian manifold $(M,g)$, $n \geq 4$,
satisfies on ${\mathcal{U}}_{S} \cap {\mathcal{U}}_{C} \subset M$
(\ref{pseudo}) and (\ref{4.3.012})
or
(\ref{pseudo}) and (\ref{genpseudo01}),
then on this set we have
\begin{eqnarray}
R &=& \frac{\phi}{2}\, S\wedge S + \mu\, g\wedge S + \frac{\eta}{2}\, g \wedge g ,
\label{eq:h7a}
\end{eqnarray}
where 
$\phi$, $\mu $ and $\eta $ are some functions on ${\mathcal U}_{S} \cap {\mathcal U}_{C}$
(cf. {\cite[Section 1] {2016_DGJZ}}).
A semi-Riemannian manifold $(M,g)$, $n \geq 4$, satisfying (\ref{eq:h7a}) on 
${\mathcal U}_{S} \cap {\mathcal U}_{C} \subset M$ 
is called a {\sl Roter type manifold}, or a {\sl Roter manifold}, or a {\sl Roter space} 
\cite{{P106}, {2016_DGJZ}, {DGP-TV02}, {DHV2008}}.
If we set in Proposition 3.2 $B = R$, and in a consequence,
$(\mathrm{Ric}(B)) = S$, $(\mathrm{Ric}(B))^{2} = S^{2}$,
$\kappa (B) = \kappa$ and 
$(\mathrm{Weyl}(B)) = C$, we obtain a family of curvature conditions which are satisfied on Roter type manifolds.
In particular, (\ref{pseudo}) and (\ref{Riccipseudo07})-(\ref{genpseudo01}), as well as
\begin{eqnarray*}
C \cdot R - R \cdot C &=& Q(S,C) - \frac{\kappa }{n-1}\, Q(g,C) ,    
\end{eqnarray*} 
are satisfied on such manifolds
(see, e.g., 
{\cite[Theorem 3.2, Proposition 3.3] {2016_DGJZ}}, {\cite[Theorem 2.4] {2018_DGZ}}). 
Moreover, from (\ref{roter71}) it follows that every Roter type manifold is a partially Einstein space.  
Roter type manifolds and in particular Roter type hypersurfaces 
in semi-Riemannian spaces of constant curvature were studied in:
\cite{{DecuP-TSVer}, {P106}, {2013_DGHHY}, {2015_DGHZ}, {2011_DGP-TV}, 
{R102}, {DeKow}, {DePlaScher}, {DeScher}, {G5}, {Kow01}, {Kow2}}. 

\vspace{3mm}

\noindent
{\bf{Remark 5.1.}} (i) 
{\cite[Remark 2.5 (iv) and (v)] {2018_DGZ}}
In the standard Schwarzschild coordinates $(t; r; \theta; \phi)$, 
and the phy\-sical units ($c = G = 1$), the Reissner-Nordstr\"{o}m-de Sitter ($\Lambda > 0$), and 
Reissner-Nordstr\"{o}m-anti-de Sitter ($\Lambda < 0$) spacetimes are given by the line element (see, e.g., \cite{SH})
\begin{eqnarray}
ds^{2} &=& - h(r)\, dt^{2} + h(r)^{-1}\, dr^{2} + r^{2}\, ( d\theta^{2} + \sin ^{2}\theta \, d\phi^{2}),
\label{rns01}
\end{eqnarray}
where $h(r) = 1 - (2 M / r) + ( Q^{2} / r^{2} ) - ( \Lambda r^{3} / 3)$ 
and $M, Q, \Lambda$ are non-zero constants. The metric (\ref{rns01}) satisfies (\ref{eq:h7a}) with 
\begin{eqnarray*}
\phi &=&
\frac{3}{2}\, ( Q^{2}  -  M r)\, r^{4} \, Q^{- 4},\ \ \
\mu \ =\
\frac{1}{2}\, ( Q^{4} + 3 Q^{2} \Lambda r^{4} - 3 \Lambda M r^{5} )\, Q^{-4}, \\ 
\eta
&=&
\frac{1}{12}\, 
( 3 Q^{6} + 4 Q^{4} \Lambda r^{4} - 3 Q^{4} M r + 9 Q^{2} \Lambda ^{2} r^{8} - 9 \Lambda^{2} M r^{9} )
\, r^{-4} \, Q^{- 4} .
\end{eqnarray*}
If we set $\Lambda = 0$ in (\ref{rns01}) then we obtain the line element of the Reissner-Nordstr\"{o}m
spacetime, see, e.g., {\cite[Section 9.2] {GrifPod}} and references therein. 
It seems that the Reissner-Nordstr\"{o}m spacetime is the "oldest" example of the Roter type warped product manifold.
\newline
(ii)
{\cite[Remark 2.5 (vi)] {2018_DGZ}}
Some comments on pseudosymmetric manifolds (also called Deszcz symmetric spaces),
as well as Roter spaces, are given in {\cite[Section 1] {DecuP-TSVer}}:
"{\sl{From a geometric point of view, the Deszcz symmetric spaces may well
be considered to be the simplest Riemannian manifolds next to the real space forms.}}" 
and 
"{\sl{From an algebraic point of view, Roter spaces may well be considered to
be the simplest Riemannian manifolds next to the real space forms.}}"
For further remarks and comments we refer to \cite{LV3-Foreword}, 
as well as \cite{{DHV2008}, {HV_2007}, {HaVerSigma}, {LV1}, {LV2}}.
\newline
(iii) Recently in \cite{2018_DH} warped product manifolds, with $2$-dimensional base and with fiber of constant curvature, which are Roter
type manifolds and admit geodesic mappings were constructed. In that paper it was also stated 
that manifolds geodesically related to these warped products are also Roter type manifolds. 

%\vspace{3mm}

\section{Hypersurfaces with two distinct affine principal curvatures}

A non-degenerate hypersurface $M$ in ${\mathbb{A}}^{n+1}$, $n \geq 1$,
is called an {\sl{improper affine hypersphere}} if $\mathcal{S}$ is identically zero.
If $\mathcal{S} = \lambda \, I\!d$, where $\lambda$ is a non-zero constant, then $M$
is called a {\sl{proper affine hypersphere}}
(see, e.g., {\cite[Definition 3.3] {Nomizu-Sasaki}}).
We note that if $M$ is a proper affine hypersphere in ${\mathbb{A}}^{n+1}$, $n \geq 2$, then 
(\ref{aff03zzz}), (\ref{aff03}), (\ref{aff0404}), (\ref{D01}) and (\ref{D02}) turn into 
\begin{eqnarray*}
S\ =\ \lambda\, h,\ \ \  R^{\ast}\ =\ \frac{\lambda ^{2}}{2}\, h \wedge h,\ \ \
\mathrm{Ric}(R^{\ast})\ =\ (n-1) \lambda ^{2}\, h,\ \ \ 
R^{\ast} \cdot R^{\ast}\ =\ Q(\mathrm{Ric} (R^{\ast}), R^{\ast})\ =\ 0,
\end{eqnarray*}
respectively.

If $M$ is a non-degenerate hypersurface in ${\mathbb{A}}^{n+1}$, $n \geq 3$, 
then $R^{\ast} = \frac{1}{2}\, S \wedge S$, i.e. (\ref{aff03}), holds on $M$.
Evidently, (\ref{aff03}) is a particular case of (\ref{aff0303}).
Now, in view of Proposition 3.1 and (\ref{aff55aff55})-(\ref{aff05}), we get
\begin{theo}
On any non-degenerate hypersurface $M$ in ${\mathbb{A}}^{n+1}$, $n \geq 3$, 
the following conditions are satisfied:
$\mathrm{Ric}(R^{\ast}) = \mathrm{tr}_{h}(\mathcal{S})\, S - S^{2}$, 
i.e. (\ref{aff0404})(a); 
$R^{\ast} \cdot R^{\ast}  = Q(\mathrm{Ric} (R^{\ast}), R^{\ast})$,
i.e. (\ref{aff05}); 
$R^{\ast} \cdot \mathrm{Ric} (R^{\ast}) = Q( \mathrm{Ric}(R^{\ast}), S)$, and 
\begin{eqnarray*}
R^{\ast} \cdot \mathrm{Weyl} (R^{\ast}) &=& Q(\mathrm{Ric} (R^{\ast}), R^{\ast}) 
- \frac{1}{n-2}\, h \wedge Q( \mathrm{Ric}(R^{\ast}), S) . 
\end{eqnarray*}
\end{theo}

As it was mentioned in Section 1, 
for affine quasi-umbilical $M$  
in $\mathbb{A}^{n+1}$, $n \geq 3$,
the tensor $\mathrm{Weyl} (R^{\ast})$ vanishes \cite{1990_OV1}. 
The converse statement is true, when $n \geq 4$ \cite{1992D2}. 
\begin{theo} {\cite[Theorem 4.8] {1992D2}} 
Let $M$ be a non-degenerate hypersurface in ${\mathbb{A}}^{n+1}$, $n \geq 4$.
Then $M$ is affine quasi-umbilical if and only if the tensor $\mathrm{Weyl}(R^{\ast})$ vanishes on $M$. 
\end{theo}

In Section 5 we presented the definition of partially Einstein semi-Riemannian spaces. 
Similarly, we can introduce the definition of affine partially Einstein hypersurfaces.
A hypersurface $M$ in $\mathbb{A}^{n+1}$, $n \geq 3$, 
will be called the {\sl{affine partially Einstein hypersurface}} 
if at every point $x \in M$, 
at which the affine shape operator $\mathcal{S}_{x}$ is not proportional to the identity transformation
$I\!d_{x}$, (\ref{aff51}) is satisfied. 

As it was noted in Section 1, 
if (\ref{aff05aff}) is satisfied at a point $x \in \mathcal{U}_{\mathcal{S}} \subset M$,
i.e. $M$ is an affine quasi-umbilical at $x$, 
then (\ref{aff06aff}) holds at this point. 
Evidently, (\ref{aff06aff}) is a particular form of (\ref{aff51}).
Thus any affine quasi-umbilical hypersurface is an affine partially Einstein hypersurface. 

We assume now that (\ref{aff51}) is satisfied on a non-degenerate hypersurface $M$ in ${\mathbb{A}}^{n+1}$, $n \geq 3$.
Proposition 3.1, (\ref{aff03}), (\ref{aff05}) and (\ref{aff51}), yield
\begin{eqnarray*}
R^{\ast} \cdot R^{\ast} &=& Q(\mathrm{Ric} (R^{\ast}), R^{\ast})
\ =\ Q(  \mathrm{tr}_{h}(S)\, S - S^{2}, \frac{1}{2}\, S \wedge S)
\ =\ Q( - S^{2}, \frac{1}{2}\, S \wedge S)
\\
&=&
Q( L_{1}\, S + L\, h , \frac{1}{2}\, S \wedge S)
\ =\  L\, Q( h , \frac{1}{2}\, S \wedge S) \ =\ L\, Q( h , R^{\ast} ) .
\end{eqnarray*}
Thus we have
\begin{theo}
If the condition 
$S^{2} + L_{1}\, S + L\, h = 0$, i.e. (\ref{aff51}),
is satisfied on a non-degenerate hypersurface $M$ in ${\mathbb{A}}^{n+1}$, $n \geq 3$,
then $R^{\ast} \cdot R^{\ast} =  L\, Q( h , R^{\ast} )$, 
i.e. (\ref{aff06}),
holds on $M$.   
\end{theo}

We assume again that $M$ is an affine quasi-umbilical hypersurface in $\mathbb{A}^{n+1}$, $n \geq 3$, 
satisfying at every point $x \in M$ the condition (\ref{aff05aff}). 
Thus (\ref{aff06aff}) holds at every point of $M$. 
Note that
(\ref{aff06aff}), by (\ref{aff0404})(a) and  (\ref{aff51}), turns into
\begin{eqnarray}
\mathrm{Ric} (R^{\ast}) - \rho ( \mathrm{tr}_{h} ( \mathcal{S} ) - \rho )\, h 
&=&
(n-2)\rho \, (S - \rho \, h) .
\label{aff72aff}
\end{eqnarray}
If in addition, $\rho \neq 0$ at $x$, then (\ref{aff05aff}), (\ref{aff06aff}) and (\ref{aff72aff}) yield 
\begin{eqnarray}
\mathrm{rank}(
\mathrm{Ric} (R^{\ast}) - \rho ( \mathrm{tr}_{h} ( \mathcal{S} ) - \rho )\, h ) &=& 
(n-2) \rho \, \mathrm{rank}( S - \rho\, h) \ = \ 1 .
\label{aff82aff}
\end{eqnarray}
This,
together with Theorem 6.3, leads to the following theorem. 
\begin{theo}
Let $M$ be an affine quasi-umbilical hypersurface $M$ in $\mathbb{A}^{n+1}$, $n \geq 3$, 
and let at every point $x \in M$ the condition (\ref{aff05aff}) be satisfied. 
Then (\ref{aff06}) holds on $M$, with the function $L$ defined by
$L = L(x) = \rho ( \mathrm{tr}_{h} ( \mathcal{S} ) - (n-1) \rho )$.
Moreover, if $\rho \neq 0$ at $x$  
then
(\ref{aff82aff}) holds at this point. 
\end{theo}
Further, we have
\begin{theo}
Let $M$ be a non-degenerate hypersurface $M$ in $\mathbb{A}^{n+1}$, $n \geq 4$, 
and let $\mathcal{U}$ be the set of all points 
of ${\mathcal{U}}_{ \mathrm{Ric}(R^{\ast})} \cap {\mathcal{U}}_{ \mathrm{Weyl}(R^{\ast})} \subset M$ 
at which (\ref{aff07}) is satisfied.
\newline
(i)
If $R^{\ast} \cdot R^{\ast}\, = \, L\, Q( h , R^{\ast} )$, 
i.e. (\ref{aff06}), is satisfied on $\mathcal{U}$ then (\ref{aff08}) and 
\begin{eqnarray}
R^{\ast} &=& \frac{\phi }{2}\, ( \mathrm{Ric} (R^{\ast}) - L\, h ) \wedge  ( \mathrm{Ric} (R^{\ast}) - L\, h )  
\label{aff72}
\end{eqnarray}
hold on $\mathcal{U}$, where $\phi$ is some function on this set. 
\newline
(ii) Let the following condition be satisfied on $\mathcal{U}$
\begin{eqnarray*}
R^{\ast} &=& \frac{\phi }{2}\, \mathrm{Ric} (R^{\ast}) \wedge \mathrm{Ric} (R^{\ast}) 
+ \mu \, h \wedge  \mathrm{Ric} (R^{\ast}) 
+ \frac{\eta }{2}\, 
h \wedge h ,
%\label{ORV01}
\end{eqnarray*}
where $\phi$, $\mu $ and $\eta$ are some functions on this set. 
Then we have on $\mathcal{U}$  
\begin{eqnarray*}
R^{\ast} \cdot R^{\ast} &=& 
( (n-2) ( \mu ^{2} - \phi \eta ) - \mu ) \phi ^{-1} \, Q( h , R^{\ast} ) .
\end{eqnarray*}
\end{theo}
{\bf{Proof.}} The conditions (\ref{aff05}) and (\ref{aff06}) yield 
$Q(  \mathrm{Ric} (R^{\ast}) - L\, h , R^{\ast} ) = 0$.
From this, in view of Proposition 2.1, it follows immediately that (\ref{aff72}) holds on $\mathcal{U}$.
Furthermore, from (\ref{aff72}), in view of Proposition 3.2(iii),
we get (\ref{aff08}).
The second assertion is an immediate consequence of Proposition 3.2(i).  
Our theorem is thus proved. 
\qed

\vspace{2mm}

The conditions (\ref{aff55aff55}), (\ref{aff52}) and Theorem 6.1 yield
\begin{theo}
Let $M$ be a non-degenerate hypersurface in ${\mathbb{A}}^{n+1}$, $n \geq 3$.
If $M$ is an affine Einstein$^{\ast}$ hypersurface then the following conditions are satisfied on $M$:
\begin{eqnarray*}
& & 
S^{2} - \frac{ \mathrm{tr}_{h}( \mathcal{S}^{2}) }{n}\, h 
\ = \  \mathrm{tr}_{h}(S)\, (S -  \frac{ \mathrm{tr}_{h}(  \mathcal{S} ) }{n}\, h ) ,
\ \ \ 
\kappa ( R^{\ast} )\, ( R^{\ast} \cdot h - Q(h,S)) \ =\ 0 ,\nonumber\\ 
& &
R^{\ast} \cdot R^{\ast} \ =\ \frac{ \kappa ( R^{\ast} ) }{n} \, Q( h , R^{\ast}) ,
\ \ \
Weyl (R^{\ast})
\ =\
\frac{1}{2}\, S \wedge S - \frac{ \kappa( R^{\ast} ) }{2 (n-1)n }\, h \wedge h .
%\label{aff53}
\end{eqnarray*}
\end{theo}

%\vspace{3mm}

\noindent
{\bf{Example 6.1.}} 
(i) Let $M$ be a 
locally strongly convex hypersurface 
in ${\mathbb{A}}^{n+1}$, $n \geq 3$.
Let $\rho_{1}, \rho_{2}, \ldots , \rho_{n}$ be 
affine principal curvatures  
at a point $x \in M$. 
Moreover, let 
$\lambda _{1} = \rho_{1} = \ldots = \rho_{k}$
and 
$\lambda _{2} = \rho_{k+1} = \ldots = \rho_{n}$, $1 \leq k < n$ and $\lambda_{1} \neq \lambda _{2}$.
Evidently, we have
\begin{eqnarray*}
\mathcal{S}_{x}^{2} - ( \lambda _{1} + \lambda _{2} )\, \mathcal{S}_{x} + \lambda _{1} \lambda _{2}\, I\!d_{x} = 0 .
\end{eqnarray*}
This, by (\ref{aff03zzz}), turns into
\begin{eqnarray} 
S^{2}_{x} - ( \lambda _{1} + \lambda _{2})\, S_{x} + \lambda _{1} \lambda _{2}\, h_{x} &=& 0 .
\label{aff54}
\end{eqnarray}
Now, in view of Theorem 6.3, we have at this point
\begin{eqnarray*}
R^{\ast} \cdot R^{\ast} &=& \lambda _{1} \lambda _{2}\, Q(h, R^{\ast}) .
%\label{aff56}
\end{eqnarray*}
Further, using (\ref{aff0404}) and (\ref{aff54}) we get 
\begin{eqnarray}
\ \ \ \ \ 
( \mathrm{Ric} (R^{\ast}))_{x} &=&  \mathrm{tr}_{h}( \mathcal{S}_{x})\, S_{x} - S^{2}_{x}
\ =\ 
( (k-1) \lambda_{1} + (n-k-1) \lambda_{2})\, S_{x} + \lambda _{1} \lambda _{2}\, h_{x} ,
\label{aff57}
\end{eqnarray}
\begin{eqnarray*}
(\kappa (R^{\ast}))_{x} &=&  (\mathrm{tr}_{h}( \mathcal{S}_{x}))^{2} -  \mathrm{tr}_{h}( \mathcal{S}^{2}_{x})
\ =\ 
( k \lambda_{1} + (n-k) \lambda_{2})^{2} - k \lambda_{1}^{2} - (n-k) \lambda_{2}^{2}\nonumber\\
&=& k (k-1) \lambda_{1}^{2} + (n-k) (n-k-1) \lambda_{2}^{2} + 2 k (n-k) \lambda _{1} \lambda _{2} .
%\label{aff58}
\end{eqnarray*}

\noindent
(ii) If $k = 1$, resp., $k = n-1$, then $\mathrm{rank} (  \mathcal{S}_{x} - \lambda _{2}\, I\!d_{x})  = 1$,
resp., $\mathrm{rank}(  \mathcal{S}_{x} - \lambda _{1}\, I\!d_{x} )  = 1$. 
Thus, in both cases, $M$ is affine quasi-umbilical at $x$,
and, in a consequence, $W (R^{\ast}) = 0$ at this point.
\newline
(iii) We assume that $2 \leq k \leq n-2$ and $(k-1) \lambda_{1} + (n-k-1) \lambda_{2} \neq 0$. 
Now from Theorem 6.2 and (\ref{aff57}) it follows that the tensors
$W(R^{\ast})$ and $ \mathrm{Ric}(R^{\ast}) -  ( \kappa (R^{\ast}) / n )\, h$
are non-zero at $x$. This means that 
$x \in {\mathcal{U}}_{ \mathrm{Ric}(R^{\ast})} \cap {\mathcal{U}}_{ \mathrm{Weyl}(R^{\ast})}$.
Furthermore, (\ref{aff03}), by (\ref{aff57}), turns into
\begin{eqnarray*}
R^{\ast}_{x} &=& 
\frac{\phi _{x} }{2}\, 
( (\mathrm{Ric} (R^{\ast})){x} - \lambda _{1} \lambda _{2}\, h_{x} ) 
\wedge ( (\mathrm{Ric} (R^{\ast}))_{x} - \lambda _{1} \lambda _{2}\, h_{x} ) ,
\end{eqnarray*}
where $\phi _{x} \, =\,  ( (k-1) \lambda_{1} + (n-k-1) \lambda_{2})^{-2}$. 
\newline
(iv) Let $M$ be a locally strongly convex hypersurface in ${\mathbb{A}}^{5}$
having at every point exactly two distinct affine principal curvatures 
$\lambda _{1} = \rho_{1} = \rho_{2}$
and
$\lambda _{2} = \rho_{3} = \rho_{4}$.
A family of hypersurfaces having this property is determined in \cite{2017_S_V}.
We have 
$\mathrm{rank} (S - \lambda_{1} \, h) = \mathrm{rank} (S - \lambda_{2}\, h) = 2$,
i.e. (\ref{aff101}) holds on $M$. This means that
$M$ is an affine $2$-quasi-umbilical hypersurface.
\newline
(v) In {\cite[Main Theorem] {Cece_Li}}, among other things,
a class of locally strongly convex hypersurfaces $M$ in ${\mathbb{A}}^{n+1}$, $n \geq 3$, 
with two distinct affine principal curvatures  
$\lambda_{1} = 0$ and $\lambda_{2} \neq 0$, 
with multiplicities $1 + m_{1}$, $m_{1} \geq 1$, and $m_{2}$, respectively, 
is determined.
This class of hypersurfaces contains affine quasi-umbilical hypersurfaces, provided that $m_{2} = 1$,
as well as affine $2$-quasi-umbilical hypersurfaces, provided that $m_{2} = 2$.

\vspace{2mm}

%\noindent
Using the above presented results we can easily prove the following theorem.
\begin{theo}
Let $M$ be a 
locally strongly convex hypersurface 
in ${\mathbb{A}}^{n+1}$, $n \geq 4$.
\newline
(i) 
If at every point $x$ of $M$ the affine shape operator $\mathcal{S}_{x}$ 
has exactly two distinct
affine principal curvatures  
$\lambda _{1}$ and $\lambda _{2}$ then (\ref{aff06}) holds on $M$,
with the function $L$ defined by $L = L(x) = \lambda _{1} \lambda _{2}$.
\newline
(ii) 
If at every point $x$ of $M$ the affine shape operator $\mathcal{S}_{x}$ 
has exactly two distinct 
affine principal curvatures  
$\lambda _{1}$ and $\lambda _{2}$ of multiplicities $k$ and $n-k$, respectively, 
such that 
\begin{eqnarray} 
(k-1) \lambda_{1} + (n-k-1) \lambda_{2} \neq 0 ,\ \ \
2 \leq k \leq n - 2 , 
\label{affine000}
\end{eqnarray}
then (\ref{aff72}) holds on $M$, where $\phi$ is some function. 
\newline
(iii) 
If at every point $x$ of $M$ the affine shape operator $\mathcal{S}_{x}$ 
has exactly two distinct 
affine principal curvatures  
$\lambda _{1}$ and $\lambda _{2}$ of multiplicities $\geq 2$ such that 
$\lambda _{1} \lambda _{2} > 0$
($\lambda _{1} \lambda _{2} < 0$) or 
$\lambda _{1} \neq 0$ and $\lambda _{2} = 0$
($\lambda _{1} = 0$ and $\lambda _{2} \neq 0$)
then (\ref{aff72}) and (\ref{affine000}) hold on $M$, where $\phi$ is some function. 
\newline
(iv) If the affine shape operator $\mathcal{S}_{x}$ 
of a four dimensional hypersurface $M$ has at every point exactly two distinct 
affine principal curvatures $\lambda _{1}$ and $\lambda _{2}$, 
both of multiplicities $2$, 
then $M$ is an affine $2$-quasi-umbilical hypersurface. 
\end{theo}

\section{Hypersurfaces with three distinct affine principal curvatures}

Let a non-degenerate hypersurface $M$ in $\mathbb{A}^{n+1}$, $n \geq 4$, 
be an affine $2$-quasi-umbilical hypersurface and let 
(\ref{aff101}) be satisfied on some open non-empty subset  $\mathcal{U} \subset \mathcal{U}_{\mathcal{S}}$.
If we set 
\begin{eqnarray}
A &=& S - \rho \, h  
\label{aff103}
\end{eqnarray}
then 
(\ref{aff0404})(a)
turns into
\begin{eqnarray}
\mathrm{Ric} (R^{\ast}) &=& - A^{2} + \alpha _{1}\, A + \alpha _{2}\, h,
\label{aff104}
\end{eqnarray}
where 
\begin{eqnarray}
(a)\ \ 
\alpha _{1} \ =\ \mathrm{tr}(A) + (n-2) \rho ,
\ \ & &\ \
(b)\ \ 
\alpha _{2} \ =\ \rho ( \mathrm{tr}(A) + (n-1) \rho ) . 
\label{aff105}
\end{eqnarray}
Further, by making use of Proposition 3.4 (i), (\ref{aff104}), (\ref{aff105})
and the assumption $\mathrm{rank} (A) = 2$,
we can easily check that the following equations are satisfied on $\mathcal{U}$
\begin{eqnarray}
& &
( \mathrm{Ric} (R^{\ast}) - \alpha _{2}\, h) \wedge ( \mathrm{Ric} (R^{\ast}) - \alpha _{2}\, h) \ =\ \tau \, A \wedge A ,
\label{aff107}\\
& &
\tau \ =\ \frac{1}{2} (  (\mathrm{tr}(A))^{2} -  \mathrm{tr}(A^{2}) ) + (n-2) \rho ( \mathrm{tr}(A) + (n-2) \rho ) .
\label{aff108}
\end{eqnarray}
From (\ref{aff103}) we get $S = A + \rho \, h$.
Using this, (\ref{aff03}) and (\ref{aff101}) we obtain
\begin{eqnarray*}
A \wedge A & = &  ( S - \rho \, h ) \wedge ( S - \rho \, h ) ,\nonumber\\ 
A \wedge A & = &  S \wedge S - 2 \rho \, h \wedge S + \rho ^{2} \, h \wedge h ,\nonumber\\
S \wedge S & = &  A \wedge A + 2 \rho \, h \wedge S - \rho ^{2} \, h \wedge h ,\nonumber\\
S \wedge S & = &  A \wedge A + 2 \rho \, h \wedge ( A + \rho \, h ) - \rho ^{2} \, h \wedge h ,
\end{eqnarray*}
and
\begin{eqnarray}
R^{\ast} &=& \frac{1}{2}\, S \wedge S \ =\ \frac{1}{2}\, A \wedge A + \rho \, h \wedge A + \frac{ \rho ^{2} }{2}\, h \wedge h .
\label{aff302}
\end{eqnarray}
If the function $\tau $, defined by (\ref{aff108}),  vanishes at $x \in  \mathcal{U}$ then
$\mathrm{rank} ( \mathrm{Ric} (R^{\ast}) - \alpha _{2}\, h) \leq 1$ at this point.
If  $\tau $  is non-zero at $x  \in  \mathcal{U}$  then 
(\ref{aff107}) and (\ref{aff302}) yield
\begin{eqnarray}
R^{\ast} &=& \frac{\tau ^{-1} }{2}\, ( \mathrm{Ric} (R^{\ast}) - \alpha _{2}\, h) \wedge ( \mathrm{Ric} (R^{\ast}) - \alpha _{2}\, h)
+ \rho \, h \wedge A + \frac{ \rho ^{2} }{2}\, h \wedge h .
\label{aff303}
\end{eqnarray}
From this, by a suitable contraction and an application of (\ref{aff105})(b), we obtain
\begin{eqnarray}
\ \ \ 
(n-2) \rho \, A  &=&
(1 - \tau ^{-1} ( \kappa (R^{\ast}) - n \alpha _{2} ) ) \, ( \mathrm{Ric} (R^{\ast}) - \alpha _{2}\, h)
+ \tau ^{-1}\, ( \mathrm{Ric} (R^{\ast}) - \alpha _{2}\, h)^{2} .
\label{aff304}
\end{eqnarray}
Now (\ref{aff303}) and (\ref{aff304}) yield
\begin{eqnarray}
(n-2) \tau\, R^{\ast} &=& \frac{ n-2 }{2}\, ( \mathrm{Ric} (R^{\ast}) - \alpha _{2}\, h) \wedge ( \mathrm{Ric} (R^{\ast}) - \alpha _{2}\, h)
 + \frac{ (n-2) \tau \rho ^{2} }{2}\, h \wedge h \nonumber\\
& &
+ 
( \tau - \kappa (R^{\ast}) + n \alpha _{2})\, h \wedge ( \mathrm{Ric} (R^{\ast}) - \alpha _{2}\, h)
+ h \wedge ( \mathrm{Ric} (R^{\ast}) - \alpha _{2}\, h)^{2} .
\label{aff305}
\end{eqnarray}
We set
\begin{eqnarray*}
B &=& R^{\ast} - \frac{ \alpha _{2} }{ 2(n-1) }\,  h \wedge h . 
%\label{aff307}
\end{eqnarray*}
This by suitable contractions yields 
$\mathrm{Ric} (B) = \mathrm{Ric} (R^{\ast}) - \alpha _{2}\, h$ and
$\kappa (B) = \kappa (R^{\ast}) - n \alpha _{2}$.
Now (\ref{aff305}) takes the form
\begin{eqnarray}
(n-2) \tau \, B &=&
\frac{ n-2 }{2}\, \mathrm{Ric} (B) \wedge  \mathrm{Ric} (B) 
+ h \wedge ( \mathrm{Ric} (B))^{2}
\nonumber\\
& &
+ ( \tau - \kappa (B))\, h \wedge \mathrm{Ric} (B)
+ \frac{(n-2) \tau }{2} \left( \rho ^{2} - \frac{\alpha _{2}}{n-1} \right) h \wedge h .
\label{aff307aa}
\end{eqnarray}
From (\ref{aff307aa}), in view of Proposition 3.3(ii), we get
\begin{eqnarray}
B &=&
\frac{ \phi }{2}\, \mathrm{Ric} (B) \wedge  \mathrm{Ric} (B) 
+ \frac{ \phi }{ n-2 } \, h \wedge ( \mathrm{Ric} (B))^{2}
+ \frac{ 1 - \kappa (B) \phi }{n-2} \, h \wedge \mathrm{Ric} (B)\nonumber\\
& &
+ \frac{ ( (\kappa (B))^{2} - \mathrm{tr} (\mathrm{Ric} (B))^{2})) \phi - \kappa (B)  }{ 2 (n-2)(n-1) } \, h \wedge h ,\ \ \ 
\phi \ =\ \tau ^{-1} .
\label{aff309}
\end{eqnarray}
This, by (\ref{AAconh01}), turns into
\begin{eqnarray}
\mathrm{Weyl} (B) &=& \phi \, \left(
\frac{1 }{2}\, \mathrm{Ric} (B) \wedge  \mathrm{Ric} (B) 
+ \frac{ 1 }{ n-2 } \, h \wedge ( \mathrm{Ric} (B))^{2} \right. \nonumber\\
& &
\left.
- \frac{ \kappa (B)  }{n-2} \, h \wedge \mathrm{Ric} (B)
+ \frac{ (\kappa (B))^{2} - \mathrm{tr} ( (\mathrm{Ric} (B))^{2}) }{ 2 (n-2)(n-1) } \, h \wedge h \right) . 
\label{aff310}
\end{eqnarray}

The above presented results and Proposition 3.2 lead to the following theorem.

\begin{theo}
Let $M$ be a non-degenerate affine $2$-quasi-umbilical hypersurface in $\mathbb{A}^{n+1}$, $n \geq 4$, 
and let (\ref{aff101}) be satisfied on some open non-empty subset $\mathcal{U} \subset \mathcal{U}_{\mathcal{S}} \subset M$.
Moreover, let $\tau $ be the function defined on $\mathcal{U}$ by (\ref{aff108}). 
\newline
(i) If $\tau = 0$ at a point $x \in \mathcal{U}$ then 
$\mathrm{rank} ( \mathrm{Ric} (R^{\ast}) - \alpha _{2}\, h) \leq 1$ at this point.
\newline
(ii) If $\tau \neq 0$ at a point 
$x \in \mathcal{U}$ 
then at this point the conditions 
$\mathrm{rank} ( \mathrm{Ric} (R^{\ast}) - \alpha _{2}\, h) \geq 2$, 
(\ref{aff305}), (\ref{aff309}) and (\ref{aff310}) are satisfied.
In particular, 
if (\ref{affine003}) is satisfied at every point of $\mathcal{U}$ 
then 
(\ref{aff08})
holds on this set.
\end{theo}

\vspace{2mm}

Let $M$ be a locally strongly convex hypersurface in ${\mathbb{A}}^{n+1}$, $n \geq 4$.
If at the point $x \in M$ 
the affine shape operator $\mathcal{S}_{x}$ has 
affine principal curvatures  
$\lambda _{1}$, $\lambda _{2}$ and $\lambda _{3} = \ldots  = \lambda _{n} = \rho $ then (\ref{aff108}), 
by (\ref{aff103}), turns into 
\begin{eqnarray}
\tau &=& ( \lambda_{1} + (n-3) \rho ) ( \lambda_{2} + (n-3) \rho ) .
\label{aff108aff}
\end{eqnarray}
We assume now that $\rho = 0$.
Thus (\ref{aff108aff}) reduces to $\tau = \lambda_{1} \lambda_{2}$. 
We also have $\mathrm{rank} (S) = 2$ at $x$. Furthermore, from Proposition 3.4(i) it follows that
\begin{eqnarray*}
(S^{2} - \mathrm{tr}_{h}( \mathcal{S})\, S) \wedge (S^{2} - \mathrm{tr}_{h}( \mathcal{S})\, S )
&=& 
- \frac{1}{2} ( \mathrm{tr}_{h}( \mathcal{S}^{2}) - (\mathrm{tr}_{h}( \mathcal{S}))^{2})\, S \wedge S  
\end{eqnarray*}
at $x$. This, by (\ref{aff0404})(b), turns into
\begin{eqnarray}
\mathrm{Ric} (R^{\ast}) \wedge \mathrm{Ric} (R^{\ast})
&=& 
 \frac{1}{2} ( (\mathrm{tr}_{h}( \mathcal{S}))^{2} - \mathrm{tr}_{h}( \mathcal{S}^{2})  )\, S \wedge S . 
\label{aff91}
\end{eqnarray}
If $\mathrm{tr}_{h}( \mathcal{S}^{2}) - (\mathrm{tr}_{h}( \mathcal{S}))^{2}) \neq 0$ at given point $x$ then
(\ref{aff03}) and (\ref{aff91}) yield
\begin{eqnarray*}
R^{\ast} 
&=& ( (\mathrm{tr}_{h}( \mathcal{S}))^{2} - \mathrm{tr}_{h}( \mathcal{S}^{2}) )^{-1} \, \mathrm{Ric} (R^{\ast}) 
\wedge \mathrm{Ric} (R^{\ast}) .
%\label{aff92}
\end{eqnarray*}

The presented above facts, Proposition 3.3 and Theorem 6.5(ii) lead to the following result. 
\begin{theo}
Let $M$ be a 
locally strongly convex hypersurface 
in ${\mathbb{A}}^{n+1}$, $n \geq 4$,
and let at every point $x$ of $M$ 
the affine shape operator $\mathcal{S}_{x}$ has three distinct 
affine principal curvatures 
$\lambda _{1}$, $\lambda _{2}$ and $\lambda _{3}$ of the form 
$\lambda _{1} \neq 0$, $\lambda _{2} \neq 0$ and $\lambda _{3} = \ldots = \lambda _{n} = 0$. 
Then on $M$ we have
$R^{\ast} = \frac{\phi}{2 }\, \mathrm{Ric} (R^{\ast}) \wedge \mathrm{Ric} (R^{\ast})$,
where the function $\phi $ is defined by $\phi = \phi (x) = (\lambda _{1} \lambda _{2})^{-1}$.
Moreover, the following conditions are satisfied on $M$: $\kappa (R^{\ast}) = 2 \phi ^{-1}$,
$(\mathrm{Ric} (R^{\ast}))^{2} = ( \kappa (R^{\ast}) / 2)\, \mathrm{Ric} (R^{\ast})$,
$R^{\ast} \cdot R^{\ast} = 0$,  
$R^{\ast} \cdot \mathrm{Weyl} ( R^{\ast} ) = 0$,    
\begin{eqnarray*}
\mathrm{Weyl} ( R^{\ast} ) \cdot R^{\ast} 
&=&  
Q ( \mathrm{Ric} (R^{\ast}) - \frac{ \kappa (R^{\ast}) }{n-1}\, h, \mathrm{Weyl} ( R^{\ast} ) ),\\
\mathrm{Weyl} ( R^{\ast} ) \cdot \mathrm{Weyl} ( R^{\ast} ) 
&=&
- \frac{ (n-3) \kappa (R^{\ast}) }{2 (n-2) (n-1) }\,
Q (h, \mathrm{Weyl} ( R^{\ast} ) ) ,
\end{eqnarray*}
and (\ref{aff08}).
\end{theo}

\vspace{2mm}

Let $M$ be a non-degenerate hypersurface in $\mathbb{A}^{n+1}$, $n \geq 3$. 
From (\ref{aff0404})(a) we get 
\begin{eqnarray}
( \mathrm{Ric} (R^{\ast}) )^{2} &=& S^{4} - 2 \mathrm{tr}_{h}(\mathcal{S})\, S^{3} + ( \mathrm{tr}_{h}(\mathcal{S}))^{2}\, S^{2} . 
\label{aff701} 
\end{eqnarray}
In addition, we assume that $M$ is 
a locally strongly convex hypersurface having at every point three distinct 
affine principal curvatures  
$\lambda_{0}$, $\lambda_{1}$ and $\lambda_{2}$. Thus we have on $M$ 
\begin{eqnarray}
S^{3} &=& \alpha \, S^{2} + \beta \, S + \gamma \, h ,
\label{aff702} 
\end{eqnarray}
where $\alpha = \lambda _{0} + \lambda _{1} + \lambda _{2}$, 
$\beta =  - \lambda _{0} ( \lambda _{1} + \lambda _{2}) - \lambda _{1} \lambda _{2}$
and $\gamma = \lambda _{0} \lambda _{1} \lambda _{2}$.  
From (\ref{aff702}) we obtain
\begin{eqnarray}
S^{4} &=& ( \alpha ^{2} + \beta ) \, S^{2} + (\alpha \beta + \gamma)\, S + \alpha \gamma \, h .
\label{aff703} 
\end{eqnarray}
Now (\ref{aff701}), by making use of (\ref{aff702}) and (\ref{aff703}), turns into 
\begin{eqnarray}
A &=& \mu \, S ,
\label{aff70456} 
\end{eqnarray} 
where
\begin{eqnarray}
A &=& ( \mathrm{Ric} (R^{\ast}) )^{2} + ( ( \alpha - \mathrm{tr}_{h}(\mathcal{S}) )^{2} + \beta )\, \mathrm{Ric} (R^{\ast}) 
- \gamma (\alpha - 2 \mathrm{tr}_{h}(\mathcal{S}) )\, h ,
\label{aff707} \noindent\\
\mu &=& 
\gamma + ( \alpha - \mathrm{tr}_{h}(\mathcal{S}) )
( \beta 
+
\mathrm{tr}_{h}(\mathcal{S})  
( \alpha - \mathrm{tr}_{h}(\mathcal{S}) ) ) .
\label{aff705} 
\end{eqnarray}  
If $\mu $ vanishes at a point $x \in M$ then $M$ is 
affine partially Einstein$^{\ast}$ 
at this point.
If $\mu $ is non-zero at a point $x \in M$ then 
from (\ref{aff70456}) it follows that $S = \mu^{-1} A$ is satisfied    
on some neighbourhood $U \subset M$ of this point.  
Moreover, (\ref{aff03}) takes on $U$ the form
\begin{eqnarray} 
R^{\ast} &=& \frac{1}{2 \mu ^{2}}\, A \wedge A . 
\label{aff706} 
\end{eqnarray}

Let $M$ be a locally strongly convex hypersurface in $\mathbb{A}^{n+1}$, $n \geq 3$, 
having at every point three distinct affine principal curvatures 
$\lambda _{0}$, $\lambda _{1}$ and $\lambda _{2}$ with multiplicities $1$, $n_{1}$ and $n_{2}$, respectively. 
Now we can express (\ref{aff705}) on $M$ in the following form
\begin{eqnarray}
\ \ \ \ \
\mu &=& ( \lambda_{0} + ( n_{1} - 1 ) \lambda_{1} + ( n_{2} - 1 ) \lambda_{2} )
( 
\lambda_{1} \lambda_{2} + ( n_{1} \lambda_{1} + n_{2} \lambda_{2} )
( ( n_{1} - 1 ) \lambda_{1} + ( n_{2} - 1 ) \lambda_{2} ) 
). 
\label{aff7890}
\end{eqnarray}
In a particular case, when $n = 3$, we have 
$\alpha = \mathrm{tr}_{h}(\mathcal{S})$ and $\mu = \gamma$ on $M$. 
If at every point of $M$ all affine principal curvatures are non-zero
then (\ref{aff706}) holds on $M$. 
If at every point of $M$ at least one of its affine principal curvatures vanishes then
$M$ is affine partially Einstein$^{\ast}$ hypersurface. 
Thus we have
\begin{theo}
Let $M$ be a 
locally strongly convex hypersurface 
in ${\mathbb{A}}^{n+1}$, $n \geq 3$,
and let at every point $x$ of $M$ 
the affine shape operator $\mathcal{S}_{x}$ has three distinct 
affine principal curvatures  
$\lambda _{0}$, $\lambda _{1}$ and $\lambda _{2}$,
with multiplicities $1$, $n_{1}$ and $n_{2}$, respectively. 
Then the tensor $A$ and the function $\mu$,
defined by (\ref{aff707}) and (\ref{aff705}), respectively,
satisfy on $M$ the equations (\ref{aff70456}) and (\ref{aff7890}). Moreover,
we have:
\newline 
(i) if $\mu$ vanishes at a point $x \in M$ then $M$ is affine partially Einstein$^{\ast}$ at this point, 
\newline
(ii) if $\mu $ is non-zero at a point $x \in M$ 
then (\ref{aff706}) holds on some neighbourhood $U \subset M$ of this point.  
\end{theo} 

%\vspace{2mm}

\noindent
{\bf{Example 7.1.}} 
(i) In {\cite[Section 5] {2014_OB_LV}} a class of locally strongly convex hypersurfaces $M$ in ${\mathbb{A}}^{6}$ 
with three distinct affine principal curvatures $\lambda _{0}$, $\lambda_{1}$ and $\lambda_{2}$, 
with multiplicities $1$, $2$ and $2$, respectively, is determined.
\newline
(ii) As it was mentioned in Example 6.1(v), in {\cite[Main Theorem] {Cece_Li}}
a class of hypersurfaces with two distinct 
affine principal curvatures is determined.
Moreover, in that paper    
also a class of hyperurfaces with three non-zero distinct affine principal curvatures  
$\lambda _{0}$, $\lambda_{1}$ and $\lambda_{2}$, 
with multiplicities $1$, $n_{1}$ and $n_{2}$, respectively, is determined. 
In particular, that class of hypersurfaces contains affine $2$-quasi-umbilical hypersurfaces, 
provided that $n_{1} = 1$ and $n \geq 4$.  

\vspace{2mm}

At the end of this section we introduce notions   
of affine quasi-umbilical$^{\ast}$ and affine 2-quasi-umbilical$^{\ast}$ hypersurface.

A hypersurface $M$ in ${\mathbb{A}}^{n+1}$, $n \geq 3$, 
is called {\sl{affine quasi-umbilical}}$^{\ast}$
at a point $x \in M$ if 
\begin{eqnarray} 
\mathrm{rank} ( \mathrm{Ric} (R^{\ast}) - \rho\, h ) &=& 1 ,\ \ \mbox{for some}\ \ \rho \in \mathbb{R},  
\label{affine001z} 
\end{eqnarray}
at this point.
From (\ref{affine001z}) we get  
\begin{eqnarray*} 
( \mathrm{Ric} (R^{\ast}) - \rho\, h ) \wedge ( \mathrm{Ric} (R^{\ast}) - \rho\, h ) &=& 1 .
\end{eqnarray*}
This, by a suitable contraction, yields (\ref{affine003}). 
Thus $M$ is affine partially-Einstein$^{\ast}$ at this point.
If (\ref{affine001z}) holds at every point of $M$ then $M$ is   
an affine quasi-umbilical$^{\ast}$ hypersurface. 

A hypersurface $M$ in ${\mathbb{A}}^{n+1}$, $n \geq 3$, 
is called {\sl{affine 2-quasi-umbilical}}$^{\ast}$ at a point $x \in M$ 
if
\begin{eqnarray} 
\mathrm{rank} ( \mathrm{Ric} (R^{\ast}) - \rho\, h ) &=& 2 ,\ \ \mbox{for some}\ \ \rho \in \mathbb{R},  
\label{affine002}
\end{eqnarray}
 at this point. 
We set $A = \mathrm{Ric} (R^{\ast}) - \rho\, h$. 
Now from Proposition 3.4(i) it follows that the tensor $A$ satisfies (\ref{eqn10.1})-(\ref{eqn14.2dd}).
If (\ref{affine002}) holds at every point of $M$ then $M$ is   
an affine 2-quasi-umbilical$^{\ast}$ hypersurface. 

\section{Hypersufaces in space forms with three distinct principal curvatures}

Let $N^{n+1}(c)$, $n \geq 3$, be  
a Riemannian space of constant curvature 
$c = \frac{\widetilde{\kappa}}{n(n+1)}$, where 
$\widetilde{\kappa}$ is its scalar curvature. 
Let $M$ be a connected hypersurface isometrically immersed in $N^{n+1}(c)$.
The Gauss equation of $M$ in $N^{n+1}(c)$,  reads 
(see, e.g., \cite{{2016_DGHZhyper}, {DGP-TV02}, {DGPSS}, {2018_DGZ}, {Saw-2015}})
\begin{eqnarray}
R_{hijk} &=& H_{hk}H_{ij} - H_{hj}H_{ik} 
+ c \, G_{hijk}, \ \ \ G_{hijk} \ =\ g_{hk}g_{ij} - g_{hj}g_{ik} , 
\label{realC5}
\end{eqnarray}
where $g_{hk}$, $R_{hijk}$, $G_{hijk}$ and $H_{hk}$ 
are the local components of the metric tensor $g$, the curvature tensor $R$, the tensor $G$    
and the second fundamental tensor of $M$, respectively.
From (\ref{realC5}), by the contraction with $g^{ij}$, we get
\begin{eqnarray}
S_{hk} - (n-1) c\, g_{hk} &=& \mathrm{tr}\, (H)\, H_{hk} - H^{2}_{hk} ,
\label{realC6}
\end{eqnarray}
where 
$S_{hk}$ and $g^{hk}$ are the local components of the Ricci tensor $S$ and the tensor $g^{-1}$, respectively,  
and
$\mathrm{tr}\, (H) = g^{ij}H_{ij}$, $H^{2}_{hk} = g^{ij} H_{hi}H_{jk}$. 
Further, from (\ref{realC6}) we get immediately
\begin{eqnarray}
S^{2} - 2 (n-1) c\, S + (n-1)^{2} c^{2} g 
&=& 
H^{4} - 2 \mathrm{tr}\, (H)\, H^{3} 
+ 
(\mathrm{tr}\, (H))^{2} \, H^{2} .
\label{realC7}
\end{eqnarray}
Evidently, the hypersurface $M$ is partially Einstein if 
(\ref{partiallyEinstein}) holds at every point of ${\mathcal U}_{S} \subset M$. 

We assume now that $M$ has
at every point three distinct 
principal curvatures  
$\lambda_{0}$, $\lambda_{1}$ and $\lambda_{2}$. 
From (\ref{realC6}) it follows that ${\mathcal U}_{S} = M$. 
Further, we have on $M$ 
\begin{eqnarray}
H^{3} &=& \alpha \, H^{2} + \beta \, H + \gamma \, g ,
\label{real702} 
\end{eqnarray}
where $\alpha = \lambda _{0} + \lambda _{1} + \lambda _{2}$, 
$\beta =  - \lambda _{0} ( \lambda _{1} + \lambda _{2}) - \lambda _{1} \lambda _{2}$
and $\gamma = \lambda _{0} \lambda _{1} \lambda _{2}$.  
From (\ref{real702}) we obtain
\begin{eqnarray}
H^{4} &=& ( \alpha ^{2} + \beta ) \, H^{2} + (\alpha \beta + \gamma)\, H + \alpha \gamma \, g .
\label{real703} 
\end{eqnarray}
Now (\ref{realC7}), by making use of (\ref{realC6}), (\ref{real702}) and (\ref{real703}), turns into 
\begin{eqnarray}
A &=& \mu \, H ,
\label{real70456} 
\end{eqnarray} 
where
\begin{eqnarray}
A &=&  S^{2} + ( ( \alpha - \mathrm{tr} (H) )^{2} + \beta  - 2 (n-1) c )\, S \nonumber\\
& &
+
( 
(n-1)^{2} c^{2} 
- ( \beta + ( \alpha - \mathrm{tr} (H) )^{2} ) (n-1)c 
- \gamma (\alpha - 2 \mathrm{tr} (H) ) )\, g ,
\label{real707} \noindent\\
\mu &=& 
\gamma + ( \alpha - 
\mathrm{tr} (H) 
)
( \beta  + \mathrm{tr} (H) ( \alpha - \mathrm{tr} (H) ) ) .
\label{real705} 
\end{eqnarray}  
Thus we see that if
$\mu $ vanishes at a point $x \in M$ then $M$ is 
partially Einstein at this point.
If $\mu $ is non-zero at a point $x \in M$ then 
from
(\ref{real70456}) 
it follows that $H = \mu ^{-1} A$
is satisfied
on some neighbourhood $U \subset M$ of this point.  
Moreover, (\ref{realC5}) takes on $U$ the form
\begin{eqnarray} 
R &=& \frac{1}{2 \mu ^{2}}\, A \wedge A + c\, G . 
\label{real706} 
\end{eqnarray}
As it was stated in Section 7, 
if $\lambda _{0}$, $\lambda _{1}$ and $\lambda _{2}$ are principal curvatures with multiplicities $1$, $n_{1}$ and $n_{2}$, respectively, 
then we can express (\ref{real705}) by (\ref{aff7890}).
Thus we have
\begin{theo}
Let $M$ be a 
hypersurface 
in a space form $N^{n+1}(c)$, $n \geq 3$,
having at every point $x$ of $M$ 
three distinct 
principal curvatures  
$\lambda _{0}$, $\lambda _{1}$ and $\lambda _{2}$,
with multiplicities $1$, $n_{1}$ and $n_{2}$, respectively. 
Then the tensor $A$ and the function $\mu$,
defined by (\ref{real707}) and (\ref{real705}), respectively,
satisfy on $M$ the equations (\ref{real70456}) and (\ref{aff7890}). Moreover,
we have:
\newline 
(i) if $\mu$ vanishes at a point $x \in M$ then $M$ is partially Einstein at this point, 
\newline
(ii) if $\mu $ is non-zero at a point $x \in M$ 
then (\ref{real706}) holds on some neighbourhood $U \subset M$ of this point.  
\end{theo} 

\noindent
{\bf{Example 8.1.}} 
(i) Hypersurfaces isometrically immersed in space forms having at every point two distinct principal curvatures 
also satisfy some curvature conditions of pseudosymmetry type, 
see, e.g., {\cite[Section 5] {2016_DGHZhyper}} and {\cite[Section 3] {2018_DGZ}}. 
\newline
(ii) 
Hypersurfaces isometrically immersed in an Euclidean space 
${\mathbb{E}}^{n+1}$, $n \geq 5$, 
having at every point three distinct principal curvatures $\lambda _{0}$, $\lambda _{1}$ and $\lambda _{2}$, 
with multiplicities $1$, $n_{1}$ and $n_{2}$, respectively, were investigated among others in \cite{Saw-2015}.
Curvature properties of pseudosymmetry type of such hypersurfaces,
in the particular case when $n_{1} = n_{2} \geq 2$ were determined in {\cite[Section 4] {Saw-2015}}.
Moreover, in {\cite[Section 5] {Saw-2015}}
a class of hypersurfaces having at every point three distinct principal curvatures $\lambda _{0}$, $\lambda _{1}$ and $\lambda _{2}$, 
with multiplicities $1$, $n_{1}$ and $n_{2}$, respectively,
satisfying in addition $( n_{1} - 1 ) \lambda_{1} + ( n_{2} - 1 ) \lambda_{2} = 0$, 
i.e. 
$\lambda _{0} + \lambda _{1} + \lambda _{2} = \mathrm{tr} (H)$, 
was determined. 
Note that on those hypersurfaces (\ref{real705}) reduces to $\mu = \lambda_{0} \lambda_{1} \lambda_{2}$.  
\newline
(iii)
In {\cite[Section 4] {DGP-TV02}} curvature properties of pseudosymmetry type of a particular class 
of $2$-quasi-umbilical minimal hypersurfaces in ${\mathbb{E}}^{n+1}$, $n \geq 4$, were determined.
\newline
(iv) We refer to {\cite[Chapter 8.5.3] {M7}} and \cite{{1995_DDDVY}, {2016_DGHZhyper}, {DGPSS}, {2018_DGZ}, {G6}} 
(see also references therein)
for further results on hypersurfaces in space forms with three distinct principal curvatures 
satisfying curvature conditions of pseudosymmetry type.

\vspace{2mm}

\noindent
{\bf{Acknowledgments.}} 
The first two authors of this paper are supported 
by the Wroc\l aw University of Environmental and Life Sciences, Poland.

\vspace{5mm}

\noindent
Ryszard Deszcz and Ma\l gorzata G\l ogowska
\newline
Department of Mathematics, 
Wroc\l aw Univeristy of Environmental and Life Sciences
\newline
Grunwaldzka 53, 50-357 Wroc\l aw, Poland
\newline
{\sf{E-mail: Ryszard.Deszcz@upwr.edu.pl}},
{\sf{Malgorzata.Glogowska@upwr.edu.pl}} 

\vspace{5mm}

\noindent
Marian Hotlo\'{s}
\newline
Department of Applied Mathematics, 
Wroc\l aw Univeristy of Science and Technology
\newline
Wybrze\.{z}e Wyspia\'{n}skiego 27, 50-370 Wroc\l aw, Poland 
\newline
{\sf{E-mail: Marian.Hotlos@pwr.edu.pl}}
\end{document}